\documentclass[11pt,a4paper,twoside]{article}
\title{\bf  Large deviation estimates for exceedance times of perpetuity sequences and their dual processes}
\author{\sc Dariusz Buraczewski$^1$, Jeffrey F.\ Collamore$^\ast$, Ewa Damek, and Jacek Zienkiewicz}

\date{\it University of Wroc{\l}aw and University of Copenhagen$^\ast$}

\usepackage{latexsym}
\usepackage{amsmath}
\usepackage{amssymb}
\usepackage{amsthm}
\usepackage{chicago}
\usepackage{amsbsy}
\usepackage{multicol}
\usepackage{color}
\usepackage{graphicx}

\usepackage[utf8]{inputenc}
\usepackage{amsmath}
\usepackage{amsfonts}
\usepackage{amssymb}
\usepackage{calligra}
\usepackage{calrsfs}
\usepackage[left=2cm,right=2cm,top=2cm,bottom=2cm]{geometry}
\usepackage[mathscr]{euscript}
\DeclareMathAlphabet{\mathcalligra}{T1}{calligra}{m}{n}

\newcommand{\beq}{\begin{eqnarray}}
\newcommand{\eeq}{\end{eqnarray}}

\newcommand{\halmos}{\vspace{3mm} \hfill \mbox{$\Box$}}
\newcommand{\reals}{{\mathbb R}}
\newcommand{\E}{{\mathbb E}}

\newcommand{\integer}{{\mathbb Z}}
\newcommand{\pintegers}{{\integer_+}}

\newcommand{\dom}{{\rm dom}\:}

\newcommand{\bfdelta}{{\boldsymbol \delta}}
\newcommand{\bfepsilon}{{\boldsymbol \varepsilon}}

\newcommand{\supp}{\mathrm{supp}}
\newcommand{\8}{\infty}
\renewcommand{\a}{\alpha}
\renewcommand{\b}{\beta}
\renewcommand{\d}{\delta}
\newcommand{\s}{\sigma}
\newcommand{\g}{\gamma}
\newcommand{\eps}{\epsilon}

\renewcommand{\P}{{\mathbb P}}

\newtheorem{thm}{Theorem}[section]

\newtheorem{Lm}{Lemma}[section]
\newtheorem{lem}{Lemma}[section]

\newtheorem{Rk}{Remark}[section]
\newtheorem*{Assert}{Assertion}

\newcommand{\Pf}{\noindent {\bf Proof.}\:\:}

\begin{document}
\maketitle

\footnotetext[1]{Corresponding author.
\par \hspace*{.1cm} All authors were supported by the NCN Grant UMO-2011/01/M/ST1/04604.
\par \hspace*{.1cm} {\it AMS 2010 subject classifications.  Primary 60H25; secondary 60K05, 60F10, 60J10, 60G70, 60K25, 60K35. }
\par  \hspace*{.1cm} {\it Keywords and phrases.  Random recurrence equations, stochastic fixed point equations, first passage times, Siegmund duality, asymptotic behavior, ruin probabilities.}}

\begin{abstract}
In a variety of problems in pure and applied probability, it is of relevant to study the large exceedance probabilities of the perpetuity sequence $Y_n := B_1 + A_1 B_2 + \cdots + (A_1 \cdots A_{n-1}) B_n$, where $(A_i,B_i) \subset (0,\infty) \times \reals$.  Estimates for the stationary tail distribution
of $\{ Y_n \}$ have been developed in the seminal papers of Kesten (1973) and Goldie (1991).  Specifically, it is well-known that if $M := \sup_n Y_n$,
then ${\mathbb P} \left\{ M > u \right\} \sim {\cal C}_M u^{-\xi}$ as $u \to \infty$.
While much attention has been focused on extending this estimate, and related estimates, to more general processes, little work has been devoted to understanding the path behavior of these processes.
In this paper, we derive sharp asymptotic estimates for the large exceedance times of $\{ Y_n \}$.  Letting $T_u := (\log\, u)^{-1} \inf\{n:  Y_n > u \}$ denote the normalized first passage time,
we  study ${\mathbb P} \left\{ T_u \in G \right\}$ as $u \to \infty$ for sets $G \subset [0,\infty)$.   We show, first, that the scaled sequence $\{ T_u \}$  converges
in probability to a certain constant $\rho > 0$.  Moreover, if $G \cap [0,\rho] \not= \emptyset$, then
${\mathbb P} \left\{ T_u \in G \right\} u^{I(G)} \to C(G)$ as $u \to \infty$ for some ``rate function'' $I$ and constant $C(G)$.
On the other hand, if $G \cap [0,\rho] = \emptyset$, then we show that the tail behavior is actually quite complex, and different asymptotic regimes are possible.
We conclude by extending our results to the corresponding forward process, understood in the sense of Letac (1986), namely, the reflected process $M_n^\ast := \max\{ A_n M_{n-1}^\ast + B_n, 0 \}$
for $n \in \pintegers$, where $M_0^\ast=0$.
Using Siegmund duality, we relate the first passage times of $\{ Y_n \}$ to the finite-time exceedance probabilities of  $\{ M_n^\ast \}$, yielding a new result concerning the convergence
of $\{ M_n^\ast \}$ to its stationary distribution.
\end{abstract}

\section{Introduction}

Since the pioneering work of \shortciteN{HK73} and \shortciteN{WV79}, there has been continued interest in the probabilistic study of perpetuity sequences.
Much of this interest has been driven by a wide variety of applications.  Perpetuity sequences arise naturally
in connection with the ARCH and GARCH financial time series models (\shortciteN{RE82}, \shortciteN{TB86}, \shortciteN{TM03}), the
Asian options in discrete and continuous time (\shortciteN{HGMY93}, \shortciteN{PCFPMY01}),
and in insurance mathematics (\shortciteN{JP02},  \shortciteN{CKRK08}, \shortciteN{JC09}).  From a theoretical perspective, they appear
naturally in connection with the weighted branching process and branching random walk (\shortciteN{AI2009}, \shortciteN{bur2009}, \shortciteN{gui1990}, \shortciteN{QL00}).   Indeed, utilizing an
argument in \shortciteN{gui1990} and \shortciteN{QL00}, it is possible to relate the tail behavior of a perpetuity sequence to that of an associated nonhomogeneous recursion,
leading to further applications, for example, to Mandelbrot cascades (\shortciteN{BDGM2014}) and the Quicksort algorithm in computer science.

A central issue arising in all of these problems is the characterization of the tail behavior of the perpetuity sequence.  Namely,
letting $\{ (A_i, B_i):  i \in \pintegers \}$ be an i.i.d.\ sequence of random variables taking values in $(0,\infty) \times \reals$,  and letting
\begin{equation} \label{intro1}
Y_n = B_1 + A_1 B_2 + \cdots + (A_1 \cdots A_{n-1}) B_n, \quad n=1,2,\ldots,
\end{equation}
then it is of interest to consider
\begin{equation} \label{intro2}
{\mathbb P} \left\{ V > u \right\} \quad \mbox{as} \:\: u \to \infty,
\end{equation}
where, typically,
\[
V:= \lim_{n \to \infty}Y_n \quad \mbox{or} \quad V := \sup_n Y_n.
\]
In either case, it is well-known that under mild regularity conditions,
\begin{equation} \label{goldie-intro}
{\mathbb P} \left\{ V > u \right\} \sim {\cal C} u^{-\xi} \quad \mbox{as} \quad u \to \infty
\end{equation}
for certain positive constants ${\cal C}$ and $\xi$ (cf.\ \shortciteN{HK73}, \shortciteN{CG91}).

Much recent work has been devoted to showing that the estimate in \eqref{goldie-intro} extends well beyond the setting of perpetuity sequences.
Following \shortciteN{GL86}, it is helpful to first observe that $\{ Y_n \}$ can be identified as the {\it backward process} generated
by the affine map $\Phi(x) = Ax + B$, where $(A,B) \stackrel{\cal D}{=} (A_1,B_1)$ (and $\stackrel{\cal D}{=}$ denotes equality in distribution).
More precisely, letting $\Phi_i(x) = A_i x + B_i$ for $i \in \pintegers,$ then
\begin{equation} \label{intro3}
Y_n = \Phi_1 \circ \cdots \circ \Phi_n(0), \quad n=1,2,\ldots.
\end{equation}
The limiting behavior of this sequence is, of course, the same as that of the corresponding {\it forward process}, namely,
\begin{equation} \label{intro4}
 Y_n^* := \Phi_n \circ \cdots \circ \Phi_1(0), \quad n =1,2,\ldots.
\end{equation}
Then it is natural to consider more general random functions, including Markov-dependent sequences and random matrices.
 Extensions
of this type can be also found, for example, in recent work of \shortciteN{GA03}, \shortciteN{AM2012}, \shortciteN{BB2014},
\shortciteN{DarEwa2009},
\shortciteN{JC09}, Collamore and Vidyashankar (2013a,b), \shortciteN{GL2013},
\shortciteN{CKSP04},
\shortciteN{NECSOZ09},
\shortciteN{MM11},   \shortciteN{AR07}.
We note that for the process \eqref{intro4}, recursions generated by random matrices were also considered in Kesten's (1973) original work.
Moreover, some refined large deviation asymptotics for related recursive structures can be found in \shortciteN{BDMZ} and \shortciteN{BDZ}.

In contrast, very little is known concerning the path properties of perpetuity sequences.
Two natural questions, well motivated by the theory of random walks, are the characterization of
the distribution of the first passage time of the sequence in \eqref{intro3}, and the convergence of the sequence in \eqref{intro4} to its stationary distribution.  Indeed, these two questions are very much
the same, since it is known by extensions of classical duality for random walks that
\begin{equation} \label{intro5}
{\mathbb P} \left\{ Y_k > u, \mbox{ some } k \le n \right\} = {\mathbb P} \big\{ { M_n^*} > u \big\},
\end{equation}
where $\{ { M_n^*} \}$ is defined as in \eqref{intro4}, but with $\Phi(x)$ replaced with $\widetilde{\Phi}(x) := (Ax+B)^+$, and $Y_0=0= M^*_0$.
(Cf.\ \shortciteN{DS76a}, \shortciteN{SAKS96}, and the discussion in Section \ref{section: 2} below.)
Thus, the {\it finite}-time exceedances of $\{ M^*_n \}$ can be analyzed
through the first passage times of $\{ Y_n \}$, and vice versa.

The primary objective of this article is to study the asymptotic distribution of the scaled first passage time
\[
T_u := \frac{1}{\log u} \inf \left\{ n:  Y_n > u \right\} \quad \mbox{as} \:\: u \to \infty.
\]
Motivated by the large deviation theory for random walks,  developed in the classic papers of
Donsker and Varadhan (cf.\ \shortciteN{SV84}), 
we study the asymptotic behavior of
\[
{\mathbb P} \left\{ T_u \in G \right\} \quad \mbox{as} \:\: u \to \infty, \quad \mbox{where} \:\: G \subset [0,\infty).
\]
We begin by showing that, conditional on $\{ T_u < \infty \}$,
\[
T_u \to \rho \quad \mbox{in probability}
\]
for some positive constant $\rho$,
thus describing the ``most likely'' first passage time into the set $(u, \infty)$.
We then characterize the asymptotic distribution of $\{T_u\}$ on the respective time intervals, $[0,\rho]$ and $[\rho,\infty)$, where the analysis on these
two regions turns out to be quite different.   On the first of these regions, namely $[0,\rho]$, we show that
there exists a ``rate function'' $I:  [0,\infty) \to [0,\infty)$ such that
 \begin{equation} \label{intro9}
{\mathbb P} \big\{ T_u \le \tau \big\} \sim  \left\{ \begin{array}{l@{\quad, \quad}l}  \left( C(\tau)/\sqrt{\log u} \right) \: u^{-I(\tau)} & \tau < \rho,\\[.4cm]
\left( {\cal C}/2 \right) u^{-\xi} & \tau = \rho,\\[.4cm]
{\cal C} u^{-\xi} & \tau > \rho,
\end{array} \right.
\end{equation}
where $\{ C(\tau):  \tau \in \reals\}$ is a collection of positive constants and $({\cal C},\xi)$ is given as in \eqref{goldie-intro}.  (See Theorems \ref{thm:21} and \ref{thm:22} below.)
In \eqref{intro9}, the case $\tau < \rho$ describes the ``small-time" behavior of $\{ T_u \}$, while the case $\tau = \rho$ can be viewed as the ``critical case,"
which, from a mathematical perspective, requires a much more elaborate analysis.  We note that
\eqref{intro9} is a considerable  refinement of \shortciteN{HN01}, who gave initial estimates for these probabilities in the ``small-time" case, namely,
rough logarithmic asymptotics for
 $\log  {\mathbb P} \big\{ T_u \le \tau \big\}$ as $u \to \infty$ when
$\tau < \rho$.  

As with the critical case, the asymptotic distribution of $\{ T_u \}$ for large times, when $\tau > \rho$, is also complex, requiring new mathematical techniques.
Indeed, as we demonstrate, these asymptotics can be quite different from those expected from the large deviation theory of random walks, which, based on
 \shortciteN{GA55}, Asmussen (2000, Chapter 4), and \shortciteN{JC98}, would suggest that
\begin{equation} \label{intro10}
{\mathbb P} \big\{\tau \le T_u < \infty \big\} \sim \frac{C(\tau)}{\sqrt{\log u}} u^{-I(\tau)} \quad \mbox{as} \quad u \to \infty.
\end{equation}
As we  show, under certain  conditions,  the previous formula fails to hold and we obtain very different asymptotic behavior, not only for ${\mathbb P} \big\{\tau \le T_u < \infty \big\}$, but also for $\log {\mathbb P} \big\{\tau \le T_u < \infty \big\}$; thus, even the polynomial decay rate suggested by \eqref{intro10} need not hold, in general.
Indeed, in Theorems \ref{thml1} and \ref{thml2} below, we provide asymptotic estimates showing that under certain conditions,
\begin{equation} \label{new-intro.1}
\limsup _{u\to \8}\log {\mathbb P} \big\{\tau \le T_u < \infty \big\}\leq -I(\tau ),
\end{equation}
while under other assumptions,
\begin{equation} \label{new-intro.2}
\liminf _{u\to \8}\log {\mathbb P} \big\{\tau \le T_u < \infty \big\} > -I(\tau ).
\end{equation}
In this way, we exhibit an interesting asymmetry between the large-time behavior and the small-time behavior of $\{ T_u \}$. These last results are quite technical and
show that for $\tau >\rho $, the story is very interesting, challenging, and not fully understood.

We now turn to a more precise statement of our results.  In the process, we also connect the convergence in \eqref{intro9} to that of the dual process of $\{ Y_n \}$.
The proofs are deferred to Sections \ref{section: 3}--\ref{section: 5}, where we establish our main results, respectively, for the three asymptotic regimes ($\tau < \rho$, $\tau=\rho$,
$\tau > \rho$) which we have just described.

%
%
%
%

\setcounter{equation}{0}
\section{Statement of results}
\label{section: 2}
\subsection{A class of stochastic recursions}
Before stating our main results, we first introduce some notation related to our stochastic recursions and formulate a few of their basic properties. Let $\left\{ (A_i,B_i):  i =1,2,\ldots \right\}$
 be a sequence of i.i.d. random variables taking values in $(0,\infty) \times \reals$.  Throughout the paper, we will assume:
\begin{itemize}
\item $\E \left[ \log A \right] \in (-\infty,0)$ and $\E \left[ \log ^+|B| \right] <\8  $.
\item For every $x\in \reals $, $\P \{ Ax+B=x\} <1$, which implies, in particular, that $\P \{B =   0\} < 1$.
\end{itemize}
We will be interested in the following two processes:    the perpetuity sequence
\begin{equation}
\label{eq: perpetuity}
Y_n := B_1 + \sum_{k=2}^n A_1\ldots A_{k-1}B_k, \quad n=1,2,\ldots, \qquad Y_0 = 0,
\end{equation}
and, particularly,  the process of partial maximums of this sequence, namely,
\begin{equation}
\label{eq: max}
M_n := \max_{0\le k\le n} Y_k, \quad n=0,1,\ldots.
\end{equation}
These sequences represent the backward processes generated by the random mappings $\Phi_i(x) = A_i x+B_i$ and $\Phi_i(x) = \left( A_i x+B_i, 0 \right)^+$, respectively.  The corresponding forward processes (defined in \eqref{intro4}) are Markov chains satisfying the respective equations
\begin{align}
\label{eq: recursions}
Y_n^* &= A_n Y^*_{n-1} + B_n,\nonumber\\[.2cm]
M_n^* &= \left( A_n M_{n-1}^* + B_n \right)^+.
\end{align}
If $\E \left[ \log A \right] <0$ and
$\E \left[ \log^+ |B| \right] <\8$, then it is well-known that  $\{Y_n\}$ converges pointwise to
$$Y = \sum_{k=1}^\8 A_1\ldots A_{k-1} B_k,$$
while $M_n$ converges a.s.\ to $$M = \sup_{n\ge 0 } Y_n,$$
where $Y$ and $M$ are finite a.s.
Then $Y$ and $M$ are called stationary solutions, since they satisfy the stochastic fixed point equations
\begin{align}
\label{eq: random equation}
Y & \stackrel{\cal D}{=} AY+B,\quad Y \mbox{ independent of } (A,B);\\[.2cm]
M & \stackrel{\cal D}{=} (AM+B)^+,\quad M \mbox{ independent of } (A,B).
\end{align}

In this paper, our objective will be to describe the path behavior of $\{Y_n\}$ and $\{M_n\}$, and, in this connection, it will be of interest to compare the limiting quantities
we obtain to the tail behavior of $Y$ and $M$.    To this end, define the generating functions
\begin{align*}
\lambda(\alpha) &= {\mathbb E} \left[ A^\alpha\right], \quad
\Lambda(\alpha) = \log \lambda(\alpha), \quad \alpha \in \reals;\\[.2cm]
\lambda_B(\alpha) &= {\mathbb E} \left[ |B|^\alpha\right], \quad
\Lambda_B(\alpha) = \log \lambda_B(\alpha), \quad \alpha \in \reals.
\end{align*}
Note by the convexity of $\Lambda$ and $\Lambda_B$ that, if
$\Lambda(\a)<\8$ and $\Lambda_B(\a)<\8$ for some $\a>0$, then $\Lambda(\b)$ and $\Lambda_B(\b)$ are finite for every
$\b\in (0,\a)$.  Moreover, these functions are infinitely differentiable on the interiors of their respective domains.

We will use some fundamental properties of the solutions to the stochastic equations \eqref{eq: random equation}. First,   if $\Lambda(\a)<0$ and $\Lambda_B(\a)<\8$ for some $\alpha > 0$,
then their $\alpha$th moments must be finite, namely,
\begin{equation}
\label{eq: moments}
\E \left[ |Y|^\a \right] <\8  \quad \mbox{and} \quad \E \left[ M^\a \right] <\8;
\end{equation}
see \shortciteN{WV79}.
Next, to describe the tail behavior of $Y$ and $M$, we focus on the nonzero solution $\xi$ to the equation $\Lambda (\xi) = 0$.
More precisely, assume that for some $\xi > 0$,
\[
\Lambda(\xi) = 0, \quad \Lambda^\prime(\xi) < \infty, \quad \mbox{and} \quad \Lambda_B(\xi) < \infty.
\]
Moreover, assume that the random variable $\log \, A$ is nonarithmetic.  Then it is well-known that the tails of $Y$ and $M$ are regularly varying with index $\xi$; that is,
\begin{eqnarray}
\label{eq: asymptotic}
\P\{ Y>u \} &\sim & {\cal C}_Y u^{-\xi}\qquad \mbox{as } u\to\8; \nonumber\\[.2cm]
\P\{ M>u \} &\sim & {\cal C}_M u^{-\xi}\qquad \mbox{as } u\to\8;
\end{eqnarray}
see \shortciteN{CG91}.   Various explicit expressions for the constants ${\cal C}_Y$ and ${\cal C}_M$ are also available; see Remark \ref{rem:2.1} below.

\subsection{Main results}

Letting $\{ Y_n \}$ denote the perpetuity sequence defined in \eqref{eq: perpetuity},  and let
\begin{equation} \label{defT}
T_u := \frac{1}{\log u} \inf \left\{ n:  Y_n > u \right\}
\end{equation}
denote the scaled first passage time of $\{ Y_n \}$ into the set $(u, \infty)$.
Then our primary objective is to study the
asymptotic decay, as $u \to \infty$, of ${\mathbb P} \left\{ T_u \in G \right\}$ for $G \subset \reals$.  We will show that this probability
decays at a polynomial rate, and provide sharp asymptotic estimates describing this rate of decay.

Set
\begin{equation}
\label{eq: mu-s}
\mu(\a) = \Lambda'(\a) \quad \mbox{ and } \quad \sigma(\a) = \sqrt{\Lambda''(\a)}.
\end{equation}
To characterize the behavior of $\{ T_u \}$ as $u \to \infty$, it is helpful to first observe that, conditional on the event of ruin, the
random variable $T_u$ converges in probability to the constant $\rho = \left( \mu(\xi) \right)^{-1}$, where
$\xi$  is given as in \eqref{eq: asymptotic}.  This constant $\rho$ will play an important role in the sequel.

\begin{Lm} \label{Lm2.1}
Assume there exists $\xi > 0$ such that $\Lambda(\xi)=0$,
 and suppose that $\Lambda$ and $\Lambda_B$ are finite in a neighborhood of $\xi$ and the law of $\log A$ is nonarithmetic.
Set $\rho = \left( \mu(\xi) \right)^{-1}.$  Then for any
$\epsilon > 0$,
\begin{equation} \label{Lm2.1-1}
{\mathbb P} \left\{  \left. T_u \notin (\rho-\epsilon, \rho + \epsilon) \: \right| \: T_u < \infty \right\} \to 0 \quad \mbox{as} \quad u \to \infty.
\end{equation}
\end{Lm}
Lemma \ref{Lm2.1} will follow as a direct consequence of a stronger result, Lemma \ref{lem:4.3}, which will be proved in Section \ref{section: 4}.

Turning now to our main results, we first introduce the rate function which we will use to describe the polynomial rates of decay.  Recall that the convex conjugate
(or Fenchel-Legendre transform) of the function $\Lambda$ is defined by
$$
\Lambda^*(x) = \sup_{\a\in\reals} \{\a x-\Lambda(\a)\}, \qquad x\in \reals.
$$
Next define
$$
I(\tau) = \tau \Lambda^*\bigg(\frac 1{\tau}\bigg), \quad \tau>0,\qquad I(0)=\8.
$$
This rate function appears in the large deviation study for random walks, and is closely related to the support function in convex analysis, whose properties
are well-known (see \shortciteN{RR70}, Chapter 13).   Various convexity properties of the function $I(\cdot)$ itself (in higher dimensions) are derived in
\shortciteN{JC98}, Section 3.
Note that if we set $\tau = \left( \mu(\alpha) \right)^{-1}$ for some $\alpha \in {\rm dom}(\mu)$ (the domain of $\mu$), then it follows that
\begin{equation}
\label{eq: rate2}
I(\tau) = \a - \frac{\Lambda(\a)}{\mu(\a)};
\end{equation}
cf.\ \shortciteN{ADOZ93}, p.\ 28.

We now turn to the characterization of ${\mathbb P} \left\{ T_u \in [0,\tau] \right\}$ when $\tau < \rho$.  Recall that the function $\Lambda$ is differentiable on the
interior of its domain.  Moreover, if $\Lambda$ is also essentially smooth (namely, if we further assume that
$|\Lambda^\prime(\alpha_i)| \uparrow \infty$ for any $\{ \alpha_i \} \subset {\rm int}\, (\dom \Lambda)$ whose limit
lies on the boundary of $\dom \Lambda$), then it is well-known that $\Lambda^\prime$ maps $\reals$ onto the entire real line.
Thus, in this case, there always exists a point $\alpha(\tau)$ satisfying the equation
\begin{equation} \label{def-alpha}
\mu(\alpha(\tau)) = \frac{1}{\tau}.
\end{equation}
More generally, it is well-known that if $\tau^{-1}$ lies in the interior of the domain of $\Lambda^\ast$, then a solution $\alpha(\tau)$ exists in \eqref{def-alpha};
cf.\ \shortciteN{RE84}, Theorem VI.5.7; \shortciteN{RR70}, Theorem 23.5.

Thus, the assumption of a solution to \eqref{def-alpha} is a very weak condition, which also appears to be necessary.
In particular, when there fails to be a solution,
one usually expects to obtain only logarithmic large deviation asymptotics rather than the sharp asymptotics which are the focus of this paper.

The most important solution to \eqref{def-alpha} appears, for our purposes, when we take $\tau = \rho$, where $\rho$ is given as in the previous lemma.
Then by definition of $\rho$, we have $\alpha(\rho) = \xi$.
Then $\tau \in (0,\rho)$ if and only if $\alpha(\tau) > \xi$, which is the setting of our first main result.

\begin{thm}\label{thm:21}
Let $\tau \in (0,\rho)$ and suppose that there exists a point $\alpha \equiv \alpha(\tau) \in \reals$ such that \eqref{def-alpha} holds.
Assume that $\Lambda$ and $\Lambda_B$ are finite in a neighborhood of $\alpha$.  Then
\begin{equation} \label{thm2.1-1}
 {\mathbb P} \big\{ T_u \le \tau \big\} =  \frac{C(\tau)}{\sqrt{\log u}} u^{-I(\tau)} \left( 1 + o(1) \right) \quad \mbox{as} \quad u \to \infty,
\end{equation}
and
\begin{equation} \label{thm2.1-xxx}
 {\mathbb P} \big\{ T_u \le \tau  - L_{\tau}(u)\big\} = o\bigg( \frac{u^{-I(\tau)} }{\sqrt{\log u}}\bigg)  \quad \mbox{as} \quad u \to \infty,
\end{equation}
where $L_{\tau}(u) = \{c\log (\log u)\}/\log u$ and $c\ge \{2(\a+1)\}/\Lambda(\a)$. The constant $C(\tau)$ is given by
\begin{equation} \label{thm2.1-1a}
C(\tau) = \frac{1}{\alpha\sigma(\alpha)\sqrt{2\pi \tau}} \lim_{n \to \infty} \frac{1}{\lambda(\alpha)^n} {\mathbb E} \left[ M_n^{\alpha} \right] \in [0,\infty).
\end{equation}
Moreover, if ${\mathbb P} \left\{ A>1, B>0 \right\} > 0$, then $C(\tau)>0$.
\end{thm}
\noindent
Note that  \eqref{thm2.1-xxx} shows, heuristically, that the critical event $\{Y_n>u\}$ occurs near the end of the time interval $[0,\tau\log u]$.

Next we turn to the critical case, which arises when we take $\tau = \rho$ in the previous theorem, and compare ${\mathbb P} \left\{ T_u \le \tau \right\}$
to ${\mathbb P} \left\{ T_u < \infty \right\}$. Notice that in this case, we have $I(\rho)=\xi$, and so our rate function reduces to the decay rate described previously in \eqref{eq: asymptotic}.

\begin{thm}\label{thm:22}
Suppose that there exists $\xi >0$ such that  $\Lambda(\xi)=0$. Also, assume that $\Lambda$ and $\Lambda_B$ are finite in a neighborhood of $\xi$ and the law of $\log A$ is nonarithmetic.
Then
 \begin{equation} \label{thm2.2-1}
{\mathbb P} \big\{ T_u \le \rho \big\} = \frac{{\cal C}_M}2 \: u^{-\xi} (1+o(1)) \quad \mbox{as} \quad u \to \infty,
\end{equation}
and
 \begin{equation} \label{thm2.2-1x}
{\mathbb P} \big\{ T_u \le \rho - L_{\rho}(u)\big\} = o\big(  u^{-\xi} \big) \quad \mbox{as} \quad u \to \infty,
\end{equation}
where $L_{\rho} = b\sqrt{\{\log(\log u)\}/\log u}$ and $b> \rho \{2(\xi+1)+\rho\sigma^2(\xi)\}$,
and the constant ${\cal C}_M \in [0,\infty)$ is given as in \eqref{eq: asymptotic}.   Moreover, if
${\mathbb P}\left\{ A>1,B>0 \right\} > 0,$
then ${\cal C}_M > 0$.
\end{thm}

\begin{Rk}  {\rm
It will follow from Lemma 4.3 below that, under the conditions of the previous theorem, we also have
$$
\P\big\{ T_u \le \tau
\big\} = {\cal C}_M u^{-\xi} (1+o(1))\quad \mbox{as} \quad u\to\8, \quad \tau > \rho.
$$
}
\end{Rk}

\begin{Rk}\label{rem:2.1} {\rm  Using Goldie's (1991) original characterization, the constant ${\cal C}_M$ in Theorem 2.2 may be expressed as
\begin{equation} \label{goldie1}
{\cal C}_M = \frac{1}{\xi \mu (\xi )} {\mathbb E} \left[\left( \left(A { M  } + B \right)^+ \right)^\xi - \left( A { M} \right)^\xi \right].
 \end{equation}
Recently, certain more explicit representation formulas have been derived for ${\cal C}_M$ and ${\cal C}_Y$ in \eqref{eq: asymptotic}; see  \shortciteN{NECSOZ09}
and \shortciteN{JCAV13}.  The main representation formula in \shortciteN{JCAV13} states that, under a weak continuity assumption on $\log \,A$,
\begin{equation} \label{JCAVconst1}
{\cal C}_M = \frac{1}{\xi \mu(\xi) {\mathbb E}[\tau]} {\mathbb E}_\xi \left[ \left( V_0 + \frac{B_1}{A_1} + \frac{B_2}{A_1 A_2} + \cdots \right)^\xi {\bf 1}_{\{ \tau = \infty \}} \right],
\end{equation}
where ${\mathbb E}_\xi[\cdot ]$ denotes expectation in the $\xi$-shifted measure (defined formally in Section 3 below),
$\tau-1$ is the first regeneration time of the forward process $\{ M_n^\ast \}$ in \eqref{eq: recursions},
and $M_0^\ast$ is chosen such that $M_0^\ast \stackrel{\cal D}{=}
M^\ast_{\tau}$.  Specifically, if ${\mathbb P} \left\{ B<0 \right\} >0$, then $\tau-1$ can be taken to be the return time of $\{ V_n \}$ to the origin.
In particular, under these conditions, the positivity of ${\cal C}_M$ follows readily from \eqref{JCAVconst1}.
Moreover, under the weaker requirements of \shortciteN{JCAV13}, Theorem 2.2, together with the additional assumption
that $\{ M_n^\ast \}$ is $\psi$-irreducible (which is implicitly assumed in Section 9 of that article), one obtains \eqref{JCAVconst1} for the $k$-chain $\{ M_{kn}^\ast:  n=1,2,\ldots \}$, as well
as the alternative representation
\begin{equation} \label{JCAVconst}
{\cal C}_M = \frac{1}{\xi \mu(\xi)} \lim_{n \to \infty} \frac{1}{n} {\mathbb E} \big[ M_n^\xi \big],
\end{equation}
which is readily seen to have a closely related form to \eqref{thm2.1-1a}.}
\end{Rk}

Finally, we turn to the case where $\tau > \rho$.
Interestingly, in this case, we do {\it not} obtain a complete analog of Theorem \ref{thm:21}.  Indeed,
counterexamples can be constructed where the asymptotics differ from those one might expect from the large deviation theory for random walks, as described in \eqref{intro10}.
For $\tau > \rho$,
the condition that appears to lead to these counterexamples is that ${\mathbb E} \left[ \log A \right] > \Lambda (\alpha(\tau))$.  In this case,
the true probability may decay at a {\it slower} polynomial
rate than $I(\tau)$.   More precisely,
within a rather flexible class of processes with  ${\mathbb E} \left[ \log A \right] > \Lambda (\alpha(\tau))$, we have 
\begin{equation} \label{sec2-counter}
{\mathbb P} \left\{ \tau \le T_u < \infty \right\} \ge {\cal D}_0 u^{-I(\tau)+\delta}, \quad \text{for sufficiently large }u.
\end{equation}
On the other hand, under different hypotheses which, in particular, imply ${\mathbb E} \left[ \log A \right] <\Lambda(\a(\tau))$, we obtain that
\begin{equation}
{\mathbb P} \left\{ \tau \le T_u < \infty \right\} \le  \frac{ {\cal D}_1}{\sqrt{\log u}} \ u^{-I(\tau)}, \quad \text{for sufficiently large }u.
\end{equation}
Thus, in this case, one cannot expect a direct analog of Theorems \ref{thm:21} and  \ref{thm:22}, and our next theorem can, in effect, be viewed
as a source of counterexamples to the natural conjecture suggested by \eqref{intro10}.


\begin{thm}\label{thml1}
Let $\tau \in (\rho,\infty)$, and suppose that there exists a point $\alpha \equiv \alpha(\tau) \in {\rm int}(\dom \Lambda)$ such that \eqref{def-alpha} holds
and
\begin{equation}
\label{eq: count1}
\mu(0) = \E \left[ \log A  \right] >\Lambda(\a).
\end{equation}
Moreover, assume that $B=1$ a.s.\ and the law of $A$ has a strictly positive continuous density on ${\mathbb R}$.
Then there exist positive  constants ${\cal D}_0$ and $\delta$ such that, for sufficiently large $u$,
\begin{equation}
\P\left\{ Y_{n_u-1} \le u \mbox{ and } Y_{n_u}>u  \right\} \ge {\cal D}_0 u^{-I(\tau)+ \delta}, \qquad n_u=\lfloor \tau\log u \rfloor.
\end{equation}
\end{thm}

\begin{Rk}{\rm
Since the construction in the theorem is quite involved, we have restricted our attention to the case $B=1$; however, the theorem can also be established
under the weaker assumption that $B > 0$ a.s.  For more details, see the discussion in Section 5.1 following the proof of the theorem. }
\end{Rk}

While the previous lemma leads essentially to a negative conclusion, we also have the following complementary result.

\begin{thm}\label{thml2}
Let $\tau \in (\rho,\infty)$, and suppose that there exists a point $\alpha \equiv \alpha(\tau) \in {\rm int}(\dom \Lambda)$ such that \eqref{def-alpha} holds
and
\begin{equation}
\label{eq: count2}
\Lambda(\b) < \Lambda(\a)\quad \mbox{\rm for some } \b < \min\{1,\a\}.
\end{equation}
Assume that $B > 0$ a.s. and $\lambda_B(-\alpha)< \infty$, and assume that the law of $(A,B)$ has compact support and  that $A$ has a bounded density.
Then there exist finite constants ${\cal D}$ and $U$ such that, for all $u \ge U$, 
\begin{equation}
\P\left\{ Y_{n_u + k -1} \le u \mbox{ and } Y_{n_u+k}>u  \right\} \le \frac{{\cal D} \varrho^k}{\sqrt{\log u}}  {u^{-I(\tau)}}, \qquad n_u=\lfloor \tau\log u \rfloor,
\end{equation}
where $\varrho := \lambda(\alpha) \in (0,1)$ and $k$ is any nonnegative integer.  Thus, for sufficiently large $u$,
\begin{equation}
{\mathbb P} \left\{ \tau \le T_u < \infty \right\} \le  \frac{ {\cal D}_1}{\sqrt{\log u}} \ u^{-I(\tau)}
\end{equation}
for some positive constant ${\cal D}_1$.
\end{thm}

\begin{Rk}{\rm
In these theorems,
conditions \eqref{eq: count1} and \eqref{eq: count2} determine the relevant asymptotic regime.
At first sight, it in not immediately clear that there are  processes which satisfy these assumptions.   In fact, such processes exists quite generally; see the discussion
in Section 5 and, in particular, Lemma 5.1.
}\end{Rk}

We conclude this section by relating the previous results to the convergence of the corresponding forward sequence $\{M_n^*\}$, where $M_n^\ast := \left(A_n M_{n-1}^\ast + B_n \right)^+,$
$n=1,2,\ldots$, and $M_0^\ast = 0$.
Borrowing terminology from queuing theory, $\{ M_n^*\}$ is called the ``content process'' corresponding to the ``risk process''
\begin{equation} \label{def-U}
U_n := \left( \frac{U_{n-1}}{A_n} - \frac{B_n}{A_n} \right)^+ , \quad n=1,2,\ldots; \quad U_0 = u.
\end{equation}
Then $\{ U_n \}$ and $\{ M_n^* \}$ are dual processes in the sense of \shortciteN{DS76a}; see \shortciteN{SAKS96}, Example 6 (slightly modified). 
Following \shortciteN{SAKS96}, 
the finite-time ruin probability of $\{ U_n \}$ may be equated to the finite-time exceedance probability of $\{ M_n^\ast \}$; that is,
\begin{equation} \label{queue}
\Psi(u) := {\mathbb P} \left\{ U_k \le 0, \:\: \mbox{some } k \le n \, | \, U_0 = u \right\} = {\mathbb P} \left\{ M^*_n \ge u \right\},
\end{equation}
and a simple argument yields that $\Psi(u)$ also describes the finite-time ruin probability of $\{ Y_n \}$, namely 
$\Psi (u)= {\mathbb P} \left\{  Y_k \ge u, \:\: \mbox{some } k \le n \right\}$
(see \shortciteN{JC09}, Section 2.1).  Thus, it is natural to relate the ruin probabilities described in our previous theorems to the exceedance probabilities of $\{ M_n^\ast \}$.

In fact, the equivalence described in \eqref{queue} can be obtained more directly in our problem.
Indeed, since the finite-time distributions of the forward and backward sequences are the same, we immediately obtain that
\begin{equation}
\P\big\{ M_n^* > u\big\} =
\P\big\{ M_n > u\big\} =
\P\big\{ Y_k > u, \mbox{ some } k\le n  \big\}.
\end{equation}
Consequently,
 \begin{equation} \label{Asm-3}
 {\mathbb P} \left\{ T_u \le \tau \right\} = {\mathbb P} \left\{ M^*_{n_u} > u \right\}, \qquad n_u = \lfloor \tau \log u \rfloor.
 \end{equation}
This leads to the following theorem concerning the convergence of the process $\{ M^*_n \}$ to its stationary distribution, which now follows as an immediate
consequence of Theorems \ref{thm:21} and \ref{thm:22}.

\begin{thm}
Let $\tau \in (0,\rho)$, and suppose that there exists a point $\alpha(\tau) \in \reals$ such that \eqref{def-alpha} holds.
Assume that $\Lambda$ and $\Lambda_B$ are finite in a neighborhood of $\alpha(\tau)$.  Then for $n_u = \lfloor \tau \log u \rfloor$, we have
\begin{equation} \label{thm2.5-1}
 {\mathbb P} \big\{  M^*_{n_u}  > u \big\} =  \frac{C(\tau)}{\sqrt{\log u}} u^{-I(\tau)} \left( 1 + o(1) \right) \quad \mbox{as} \quad u \to \infty
\end{equation}
for the finite constant $C(\tau)$ in \eqref{thm2.1-1a},
and this constant is strictly positive if  ${\mathbb P} \left\{ A>1, B>0 \right\} > 0$.

Next, let $\xi$ be given as in \eqref{eq: asymptotic} and suppose that $\Lambda$ and $\Lambda_B$ are finite in a neighborhood of $\xi$ and the law of $\log A$ is nonarithmetic.
Then for $n_u = \lfloor \tau \log u \rfloor$,
 \begin{equation} \label{thm2.5-2}
{\mathbb P} \big\{ M^*_{n_u} > u \big\} =  \left\{ \begin{array}{l@{\quad, \quad}l}
\frac{1}{2} {\cal C}_M \, u^{-\xi}(1 + o(1)) & \tau = \rho,\\[.4cm]
{\cal C}_M \, u^{-\xi} (1 + o(1)) & \tau \in ( \rho,\infty],
\end{array} \right.
\end{equation}
as $u \to \infty$
for the finite constant ${\cal C}_M \in [0,\infty)$ in \eqref{eq: asymptotic}, and this constant is strictly positive if
\mbox{${\mathbb P}\left\{ A>1,B>0 \right\} > 0$.}
\end{thm}

\setcounter{equation}{0}
\section{Proof of Theorem \ref{thm:21}}
\label{section: 3}
First we introduce some further notation, as follows.  Let
 \begin{align*}
 \Pi_n &= A_1 \cdots A_n, \quad n \in \pintegers;\\
S_n &= \sum_{k=1}^n \log \, A_k =\log \Pi_n, \quad n \in \pintegers;\\
\overline{Y}_n &= \sum_{i=1}^n \Pi_{i-1} |B_i|, \quad n \in \pintegers.
\end{align*}

Also, let $\nu$ denote the probability law of $(\log A,B)$, and if $\lambda(\alpha)<\8$,
define
\begin{equation} \label{measchange}
\nu_\alpha (E) = \int_E \frac{e^{\alpha x}}{\lambda(\alpha)}d\nu(x,y), \quad E \in {\cal B}(\reals^2),
\end{equation}
where ${\cal B}(\reals^2)$ denotes the Borel sets on $\reals^2$.
Let ${\mathbb E}_\alpha[\cdot]$ denote expectation with respect to the probability measure  $\nu_\alpha$.
Note that $\mu(\alpha) := \Lambda^\prime(\alpha)$ and $\sigma^2(\alpha) := \Lambda^{\prime\prime}(\alpha)$ (defined previously in \eqref{eq: mu-s}) denote the mean
and the variance, respectively, of the random variable $\log A$ with respect to the measure $\nu_\alpha$.

We start by establishing a variant of the exponential Chebyshev inequality from large deviation theory, commonly used in conjunction with Minkowski's inequality
for perpetuity sequences
(yielding estimates which are typically not very sharp).   The next lemma will provide a sharper version of these estimates for our problem.
Before stating this result, we  recall that $\Lambda(\xi)=0$, that is, $\xi$ denotes
the critical value that determines the decay rate of ${\mathbb P} \left\{ M > u\right\}$ as $u \to \infty$.  Thus
$\lambda(\alpha) \ge 1$ for  $\alpha \ge \xi$.

\begin{Lm}
\label{lem:3.1}
Let $\alpha \ge \xi$, and assume that $\alpha$ and $\epsilon > 0$ have been chosen such that $\Lambda(\alpha + \epsilon) <\8$ and $\Lambda_B(\alpha + \epsilon) <\8$.
Then
\begin{equation} \label{lm-3.1.1}
{\mathbb P} \left\{ \overline{Y}_n > u \right\} \le \overline{C}_n \lambda^n(\alpha) u^{-(\alpha + \epsilon)}, \quad \text{for all} \:\: u >0, \quad n \in \pintegers,
\end{equation}
where
\begin{equation} \label{lm-3.1.1a}
\overline{C}_n = {b} n \left( n-1 \right)^{2(\alpha + \epsilon)} \exp\left\{(n-1) \left(\epsilon \mu(\alpha) + \epsilon^2 \sigma^2(\alpha) \right) \right\}
\end{equation}
for ${b} = \left(\pi^2/6 \right)^{\alpha + \epsilon} \left\{\lambda_B(\alpha+\epsilon)/\lambda(\alpha) \right\}$.
\end{Lm}

\Pf
From the elementary equality $\sum_{k=1}^\infty k^{-2}= \pi^2/6$, we obtain
\begin{equation}\label{lm-3.1.2}
{\mathbb P} \left\{ \overline{Y}_{n}> u \right\} \le \sum_{k=1}^n {\mathbb P} \left\{ \Pi_{k-1} |B_{k}| > \frac{6u}{\pi^2 k^2} \right\}
\le \sum_{k=1}^n {\mathbb E} \Big[ \Pi_{k-1}^{\alpha + \epsilon} |B_k|^{\alpha+\epsilon} \Big] \left( \frac{\pi^2 k^2}{6u} \right)^{\alpha + \epsilon}.
\end{equation}
Now by independence,
\[
{\mathbb E} \Big[ \Pi_{k-1}^{\alpha + \epsilon} |B_k|^{\alpha+\epsilon} \Big] = \left( {\mathbb E} \left[ A^{\alpha + \epsilon} \right] \right)^{k-1} {\mathbb E} \left[ B^{\alpha+\epsilon} \right] :=
   \left(  \lambda(\alpha+\epsilon) \right)^{k-1} \lambda_B(\alpha+\epsilon).
\]
Moreover, since the generating function $\Lambda$ is infinitely differentiable on the interior of its domain,
\[
\lambda(\alpha + \epsilon) = e^{\Lambda(\alpha + \epsilon)} \le \exp \left\{ \Lambda(\alpha) + \epsilon \mu(\alpha) + \frac{\epsilon^2{\mathfrak m}}{2}  \right\},
\]
where ${ \mathfrak m} :=  \sup \left\{ \sigma^2(\beta):  \alpha \le \beta \le \alpha + \epsilon \right\}.$   Moreover, using the continuity of the function $\sigma^2(\cdot)$, we have
that ${\mathfrak m}/2 \le \sigma^2(\alpha)$ when $\epsilon$ is sufficiently small.  Hence,
substituting the previous two equations into \eqref{lm-3.1.2}, we obtain that for sufficiently small $\epsilon$,
\begin{equation} \label{lm-3.1.3}
{\mathbb P} \left\{ \overline{Y}_{n}> u \right\} \le u^{-(\alpha + \epsilon)}  \sum_{k=1}^{n}  G(k),
\end{equation}
where
\[
G(k) = \lambda(\alpha)^{k-1} \exp\left\{ (k-1) \left( \epsilon \mu(\alpha)
+ \epsilon^2 \sigma^2(\alpha)\right) \right\} \lambda_B(\alpha + \epsilon) \left( \frac{\pi^2 k^2}{6} \right)^{\alpha+\epsilon}.
\]
 Since
$\lambda(\alpha) \ge 1$ and $\mu(\alpha):=\Lambda^\prime(\alpha) \ge 0$, it follows that $G(k)$ is increasing in $k$.
Hence $\sum_{k=1}^n G(k) \le n G(n)$, and substituting this last estimate into \eqref{lm-3.1.3} yields \eqref{lm-3.1.1}, as required.
\halmos

Next define
\[
\overline{T}_u = \frac{1}{\log u} \inf\{ n: \overline{Y}_n > u \},
\]
and note by definition that $\overline{T}_u \le T_u$ on $\{ T_u < \infty \}$.
Then as a simple consequence of the lemma, we obtain the following.

\begin{Lm}
\label{lem:3.2}
Under the assumptions of Theorem \ref{thm:21},
 \begin{equation} \label{cor-3.1.1}
 {\mathbb P} \left\{ \overline{T}_u \le \tau - L_\tau(u) \right\} = o \left( \frac{u^{-I(\tau)}}{\sqrt{\log u}} \right) \quad \text{as} \quad u \to \infty,
 \end{equation}
for any $L_\tau(u) \ge \{c \log( \log u)\}/\log u $, where $c=\left\{2(\alpha + 1) \right\}/\Lambda(\alpha).$
\end{Lm}

\Pf
Set $\zeta_u = \lfloor \log u \left( \tau - L_\tau(u)\right) \rfloor$.  Then it follows directly from the definitions that
\begin{equation} \label{cor-3.1.1a}
{\mathbb P} \left\{ \overline{T}_u \le \tau - L_\tau(u) \right\} =
{\mathbb P} \left\{ \overline{Y}_{\zeta_u} > u \right\}.
\end{equation}
Now set $\alpha \equiv \alpha(\tau)$, where $\alpha(\tau)$ is defined as in \eqref{def-alpha}.
To apply the lemma, it is helpful to first observe (using  \eqref{eq: rate2}) that
\begin{equation}  \label{cor-3.1.2a}
\left( \lambda(\alpha) \right)^{\tau \log \,u} u^{-\alpha} =  e^{ -\log u  \left({\alpha} - \tau \Lambda(\alpha) \right)} = u^{-I(\tau)}.
\end{equation}
Hence
\begin{equation} \label{cor-3.1.2}
\left( \lambda(\alpha) \right)^{\zeta_u} u^{-\alpha} \le u^{-I(\tau)} \left( \lambda(\alpha) \right)^{-L_\tau(u) \log u}.
\end{equation}
Next, choose $\epsilon \equiv \epsilon(u)$ such that  $u^{-\epsilon(u)} = (\log\, u)^{-1/2}$, which is achieved by setting
\begin{equation} \label{cor-3.1.3}
\epsilon(u) = \frac{\log \left( \sqrt{\log u} \right)}{\log u} \searrow 0, \quad u \to \infty.
\end{equation}
Then by \eqref{cor-3.1.2}, it is sufficient
to show that
\begin{equation} \label{cor-3.1.4}
\overline{C}_{\zeta_u} \left( \lambda(\alpha) \right)^{-L_\tau(u) \log u} = o(1) \quad \text{as}\quad u \to \infty,
\end{equation}
 for $\overline{C}_{\zeta_u}$ defined as in \eqref{lm-3.1.1a}.
Observe that with the choice of $\epsilon(u)$ given in \eqref{cor-3.1.3} and the upper bound $(\zeta_u -1) \le \tau \log\, u$, we obtain that
\[
\exp\left\{(\zeta_u-1) \left(\epsilon(u) \mu(\alpha) + \epsilon^2(u) \sigma^2(\alpha) \right) \right\} = O\left( \sqrt{\log u} \right) \quad \text{as} \quad u \to \infty;
\]
hence
\begin{equation} \label{cor-3.1.5}
\overline{C}_{\zeta_u} = O \left( (\log u)^{2 (\alpha+ \epsilon) + \frac{3}{2}} \right) \quad \text{as} \quad u \to \infty.
\end{equation}
Then \eqref{cor-3.1.4} follows from \eqref{cor-3.1.5}, provided that  we choose $L_\tau (u) \log u \ge c \log(\log u)$, where $c =2\left(\alpha + 1 \right)/\Lambda(\alpha).$
\halmos

From the lemma, we see that the probability of ruin in the scaled time interval $[0,\tau - L_\tau(u)]$ is negligible, so we may concentrate on the critical interval $(\tau - L_\tau(u), \tau]$.
In this region, we will argue that the process $\{ \log \,Y_n \vee 0 \}$ behaves similarly to a perturbed random walk when this process is large, that is, $\log \, Y_n$ can be approximated
by  $S_n + \bfepsilon_n$ for some perturbation term $\bfepsilon_n$ and
$S_n := \sum_{i=1}^n \log A_i$.  To analyze the behavior of the random walk $\{ S_n \}$,
the following uniform large deviation theorem, due to Petrov (1965, Theorem 2), will play a key role.

\begin{thm}[Petrov]
\label{thm:3.1Petrov} Let $a_0 = \sup_{\alpha \in {\rm dom}(\Lambda^\prime)} \Lambda'(\alpha)$.
Suppose  that $c$ satisfies $ {\mathbb E} \left[ \log A \right] < c < a_0$, and suppose that $\delta(n)$ is an arbitrary function satisfying  $\lim_{n \to \infty} \delta(n) = 0$.
Then with 
$\alpha$ chosen such that $\Lambda '(\alpha )=c$, we have that
\begin{align} \label{petrov-1}
{\mathbb P} & \left\{ S_n > n(c + \gamma_n) \right\}\nonumber\\
& \quad\quad = \frac{1}{\alpha\sigma(\alpha) \sqrt{2\pi n}} \exp\left\{-n \Big( \alpha(c+\gamma_n) - \Lambda(\alpha) + \frac{\gamma_n^2}{2\sigma^2(\alpha)}
  \left(1 + O(|\gamma_n| \right) \Big) \right\} (1+o(1))
\end{align}
as $n \to \infty$,
uniformly with respect to $c$ and $\gamma_n$ in the range
\begin{equation} \label{petrov-2}
{\mathbb E} \left[ \log \,A \right] + \epsilon \le c \le a_0 - \epsilon \quad \text{\rm and} \quad |\gamma_n| \le \delta(n),
\end{equation}
where $\epsilon >0$.
\end{thm}

\begin{Rk}{\rm
In \eqref{petrov-2}, we may have that $\sup \{ \alpha : \alpha \in \dom (\Lambda) \} = \infty$ or ${\mathbb E} \left[ \log \,A \right] =-\infty$.  In these cases, the quantities
$\infty -\epsilon$ or  $-\infty -\epsilon$ should be interpreted as arbitrary positive, respectively negative, constants.}
\end{Rk}

\noindent
{\bf Proof of Theorem \ref{thm:21}.}
\textbf{\textit{Step 1.}}  Eq.\ \eqref{thm2.1-xxx} was established in Lemma \ref{lem:3.2}; thus,
it is sufficient to show that
\begin{equation} \label{pf-th2.1.2}
{\mathbb P} \left\{  \tau - L_\tau(u) < T_u \le \tau \right\}=  \frac{C(\tau)}{\sqrt{\log u}} u^{-I(\tau)} \left( 1 + o(1) \right) \quad \mbox{as} \quad u \to \infty,
\end{equation}
for $L_\tau(u) = \left\{c \log( \log u )\right\}/\log u $, where $c=\left\{2\left(\alpha + 1 \right)\right\}/\Lambda(\alpha)$.
Indeed, by Lemma \ref{lem:3.2},
\begin{equation} \label{pf-th2.1.1}
 {\mathbb P} \left\{ T_u \le \tau - L_\tau(u) \right\} = o \left( \frac{u^{-I(\tau)}}{\sqrt{\log u}} \right) \quad \text{as} \quad u \to \infty.
 \end{equation}

Set
\[
\zeta_u = \lfloor \log u \left( \tau - L_\tau(u)\right) \rfloor
\quad \text{and} \quad \tau_u =  \lfloor \tau \log u \rfloor,
\]
and define
\[
{\cal M}_{u} = \max_{\zeta_u  < n \le \tau_u} \Big\{ B_{\zeta_u+1} + A_{\zeta_u+1} B_{\zeta_u+2} + \cdots
  + \left( A_{\zeta_u+1} \cdots A_{n-1} \right) B_n \Big\} \vee 0.
\]
Then on $\left\{ \omega \in \Omega:  \max_{\zeta_u < n \le \tau_u} Y_n(\omega) > Y_{\zeta_u}(\omega) \right\}$, we have
\begin{equation}\label{pf-th2.1.3}
\max_{\zeta_u < n \le \tau_u} Y_n  = Y_{\zeta_u} + \Pi_{\zeta_u} {\cal M}_{u},
\end{equation}
and our objective is to show that ${\mathbb P} \left\{ \max_{\zeta_u < n \le \tau_u} Y_n > u \right\}$ decays at the rate specified on the right-hand side of \eqref{pf-th2.1.2}.

 \textbf{\textit{Step 1a.}}
We begin by analyzing the second term of the right-hand side of \eqref{pf-th2.1.3}.  Observe that
\begin{equation} \label{pf-th2.1.4}
{\mathbb P} \left\{ \Pi_{\zeta_u} {\cal M}_u > u \right\} = \int_{\reals} {\mathbb P} \left\{\log  \Pi_{\zeta_u} > \log u - s \right\} dF_u(s),
\end{equation}
where $F_u$ denotes the probability distribution function of $\log {\cal M}_u$.  To evaluate this integral,
note that $\log \, \Pi_{\zeta_u} := \sum_{k=1}^{\zeta_u} \log \, A_i := S_{\zeta_u},$
and thus
\begin{equation} \label{new-P}
\log \Pi_{\zeta_u} > \log u - s \Longleftrightarrow \frac{S_{\zeta_u}}{\zeta_u} >  \frac{\log\,u - s}{\zeta_u} =: \frac{1}{\tau} + \gamma_u.
\end{equation}
Letting $\gamma_n$ be defined as in this last equation and utilizing the definition of $\zeta_u$, we then obtain
\begin{equation} \label{pf-th2.1.4a}
\zeta_u \gamma_u = \frac{L_\tau(u)}{\tau} \log u - s + \bfdelta_u, \quad \mbox{where} \:\: |\bfdelta_u| \le \frac{1}{\tau}.
\end{equation}
Consequently,
\begin{equation} \label{pf-th2.1.4b}
\gamma_u = \frac{1}{\zeta_u} \left( \frac{L_\tau(u)}{\tau} \log u - s + \bfdelta_u \right)\quad \mbox{and}\quad
\zeta_u \gamma_u^2 = \frac{1}{\zeta_u} \left( \frac{L_\tau(u)}{\tau} \log u - s + \bfdelta_u \right)^2.
\end{equation}
From these equations, it is apparent that $\gamma_u \to 0$ and $\zeta_u \gamma_u^2 \to 0$ as $u \to \infty$ and, moreover,
this convergence is uniform in $s$ provided that
$s \in [-(\log u)^{1/3},(\log \,u)^{1/3}]$.

Now set $\alpha \equiv \alpha(\tau)$ for the remainder of the proof.
Then by applying Theorem \ref{thm:3.1Petrov}, we obtain that
\begin{equation} \label{pf-th2.1.5}
{\mathbb P} \left\{ \log \Pi_{\zeta_u} > \log u - s \right\} = \frac{1}{\alpha\sigma(\alpha) \sqrt{2\pi \tau \log u}} u^{-\alpha}e^{\alpha s} \left( \lambda(\alpha) \right)^{\zeta_u} (1 + o(1)) \quad \text{as}\quad u \to \infty,
\end{equation}
uniformly in $s$ such that $\log\,s \in [-(\log \, u)^{1/3},(\log \, u)^{1/3}]$.  Letting ${\cal G}_u = \big\{\omega \in \Omega:  \log {\cal M}_u(\omega) \in  [-(\log \, u)^{1/3},$ $(\log \, u)^{1/3}] \big\}$ and returning to \eqref{pf-th2.1.4}, we then obtain
\begin{equation} \label{pf-th2.1.6}
{\mathbb P} \left\{ \Pi_{\zeta_u} {\cal M}_u > u, \: {\cal G}_u \right\}
  = \frac{1}{\alpha\sigma(\alpha) \sqrt{2\pi \tau \log u}}  \left( \lambda(\alpha) \right)^{\zeta_u} u^{-\alpha}
  {\mathbb E} \left[ {\cal M}_u^\alpha {\bf 1}_{{\cal G}_u} \right] (1 + o(1)) \quad \text{as}\quad u \to \infty.
\end{equation}
Now recall (cf.\ \eqref{cor-3.1.2a}) that
\[
\left( \lambda(\alpha) \right)^{\tau \log\,u} u^{-\alpha} = u^{-I(\tau)}.
\]
Moreover, since ${\cal M}_u \stackrel{\cal D}{=} \max \left\{Y_i:  0 \le i \le \tau_u - \zeta_u \right\} \equiv M_{\lfloor \tau_u - \zeta_u \rfloor}$, we have
\[
\lim_{u \to \infty} \frac{1}{(\lambda(\alpha))^{\tau_u-\zeta_u}} {\mathbb E} \left[ {\cal M}_u^\alpha {\bf 1}_{{\cal G}_u} \right] = \lim_{n \to \infty} \frac{1}{\lambda^n(\alpha)} {\mathbb E} \left[ M_n^\a
 {\bf 1}_{{\cal H}_n} \right],
\]
where ${\cal H}_n :=  \left\{\omega \in \Omega:  \log \, \left( M_n(\omega) \right) \in  [-e^{n/3c}, e^{n/3c} ] \right\}$.  [In the definition of ${\cal H}_n$, we have used that
$\tau_u - \zeta_u \sim L_\tau(u) \log\:u = c \log(\log u)$.]  Substituting these last two equations into \eqref{pf-th2.1.6} yields
\begin{equation} \label{pf-th2.1.7}
{\mathbb P} \left\{ \Pi_{\zeta_u} {\cal M}_u > u, \: {\cal G}_u \right\}
  = \frac{\hat{C}(\tau)}{\sqrt{\log\,u}} u^{-I(\tau)} \left(1 + o(1) \right),
\end{equation}
where
\begin{equation} \label{pf-th2.1.7a}
\hat{C}(\tau) = \frac{1}{\alpha\sigma(\alpha) \sqrt{2\pi \tau}}  \lim_{n \to \infty} \frac{1}{\lambda^n(\alpha)} {\mathbb E} \left[ M_n^\a {\bf 1}_{{\cal H}_n} \right].
\end{equation}

To complete the proof, we now show that the restriction to the sets ${\cal G}_u$ and ${\cal H}_n$ can be removed on the left- and
right-hand sides of \eqref{pf-th2.1.7}, \eqref{pf-th2.1.7a}, and that the limit in $n$ on the right-hand side of \eqref{pf-th2.1.7a} exists and is both positive and finite.
To this end, first
observe by Chebyshev's inequality that
\begin{align} \label{pf-th2.1.8}
{\mathbb P} &\left\{ \Pi_{\zeta_u} {\cal M}_u > u, \:\: \log  {\cal M}_u < -(\log u)^{1/3} \right\} \le {\mathbb P} \left\{ S_{\zeta_u} > \log u + (\log\, u)^{1/3} \right\} \nonumber\\
& \hspace*{1cm} \le \exp\left\{ - \alpha \left( \log u + (\log u)^{1/3} \right) \right\} \left( \lambda(\alpha) \right)^{\zeta_u} = o \left(\frac{1}{\sqrt{\log u}} \left( \lambda(\alpha) \right)^{\zeta_u} u^{-\alpha} \right),
\end{align}
since $\lim_{u \to \infty} \sqrt{\log  u}\: \exp\left\{-\alpha (\log u)^{1/3}\right\}= 0$.  This shows that the restriction
to values $\big\{\omega \in \Omega:  \log {\cal M}_u(\omega) \ge$ $ -(\log u)^{1/3} \big\}$ can now be removed on the left-hand side of \eqref{pf-th2.1.6}, hence the left-hand side
of \eqref{pf-th2.1.7}.

Moreover, repeating the argument leading to \eqref{pf-th2.1.7}, we find that
${\mathbb P} \left\{ \Pi_{\zeta_u} {\cal M}_u > u, \: \log {\cal M}_u > (\log u)^{1/3} \right\}$ is equal to the right-hand side of \eqref{pf-th2.1.7}, but with
${\mathbb E} \left[ M_n^\a {\bf 1}_{{\cal H}_n} \right]$ replaced with
\[
{\mathbb E} \left[ M_n^\a {\bf 1}_{{\cal H}_n^\prime} \right], \quad\text{where} \quad {\cal H}_n^\prime
:=  \left\{\omega \in \Omega:  \log \,  M_n(\omega)  > e^{n/3c}\right\}.
\]
We claim that
\begin{equation} \label{pf-th2.1.9}
\lim_{n \to \infty} \frac{1}{ \lambda^n(\alpha) } {\mathbb E} \left[ M_n^\a {\bf 1}_{{\cal H}_n^\prime} \right] = 0.
\end{equation}

Set ${\cal H}^\prime_{n,k} = \left\{\omega \in \Omega:  \log \, M_n(\omega) -e^{n/3c} \in (k-1,k]\right\}$, $k=1,2,\ldots;$
thus $\bigcup_{k \in \pintegers} {\cal H}_{n,k} = {\cal H}_n^\prime$.  Then apply Lemma \ref{lem:3.1} to obtain that
\begin{eqnarray*}
\frac{1}{\lambda^n(\alpha)} \sum_{k = 1}^\infty {\mathbb E} \left[ M_n^\a {\bf 1}_{{\cal H}^\prime_{n,k}} \right]
&{ \le }& \frac{1}{\lambda^n(\alpha)} \sum_{k =1}^\infty e^{\alpha k} \exp\big( \alpha e^{n/3c}\big)  {\mathbb P} \left\{ \overline{Y}_n > e^k \exp\big(e^{n/3c} \big) \right\}\\
&=& \overline{C}_n \exp\big( -\epsilon e^{n/3c} \big) \sum_{k=1}^\infty e^{-\epsilon k}
\end{eqnarray*}  for $\epsilon>0$ sufficiently small.
Now choose $\epsilon \equiv \epsilon(n) = n^{-2}$.  With this choice of $\epsilon(n)$, note that $\overline{C}_n = O(n^{2\alpha+1})$ and $\sum_{k=1}^\infty e^{-\epsilon k}=O(n^2)$.
Then $\overline{C}_n \exp\big( -n^{-2}e^{n/3c}\big)n^2\to 0$ as
$n \to \infty$. Thus
we obtain \eqref{pf-th2.1.9}.

From \eqref{pf-th2.1.9}, we conclude { that the restrictions } on large values can be removed in
\eqref{pf-th2.1.7} and \eqref{pf-th2.1.7a}  (that is,  the restriction that $\log {\cal M}_u \le (\log u)^{1/3}$ in \eqref{pf-th2.1.7},
and the restriction that $M_n \le e^{n/3c}$ in \eqref{pf-th2.1.7a}).
Moreover, by a trivial calculation,  the restriction to values $M_n \ge -e^{n/3c}$ can also be removed in \eqref{pf-th2.1.7a}.
Consequently, we conclude that \eqref{pf-th2.1.7} and \eqref{pf-th2.1.7a} hold {\it without} including the term ${\cal G}_u$ in \eqref{pf-th2.1.7}, or the term
${\bf 1}_{{\cal H}_n}$ in \eqref{pf-th2.1.7a}.

 \textbf{\textit{Step 1b.}}
Finally, to establish \eqref{thm2.1-1}, recall that
$\max_{\zeta_u < n \le \tau_u} Y_n  = Y_{\zeta_u} + \Pi_{\zeta_u} M_{u}$; cf.\ \eqref{pf-th2.1.3}.  Now we have just shown that
\begin{equation} \label{new-th2.1.2}
{\mathbb P} \left\{ \Pi_{\zeta_u} {\cal M}_u > u \right\}
  = \frac{\hat{C}(\tau)}{\sqrt{\log\,u}} u^{-I(\tau)} \left(1 + o(1) \right),
\end{equation}
where
\begin{equation} \label{pf-th2.1.10a}
C(\tau) = \frac{1}{\alpha\sigma(\alpha) \sqrt{2\pi \tau}}  \lim_{n \to \infty} \frac{1}{\lambda^n(\alpha)} {\mathbb E} \left[ M_n^\a \right].
\end{equation}
Moreover,
 by another application of Lemma \ref{lem:3.2}, we have that
\begin{equation} \label{new-th2.1.1}
{\mathbb P} \left\{ |Y_{\zeta_u}| > u \right\} =  o \left( \frac{u^{-I(\tau)}}{\sqrt{\log u}} \right) \quad \text{as} \quad u \to \infty.
\end{equation}
Note that \eqref{new-th2.1.1} implies the existence of a function $\Delta(u) \downarrow 0$ such that
\begin{equation} \label{new-th2.1.3}
{\mathbb P} \left\{ |Y_{\zeta_u}| > \Delta(u)u \right\} =  o \left( \frac{u^{-I(\tau)}}{\sqrt{\log u}} \right) \quad \text{as} \quad u \to \infty.
\end{equation}
Moreover, on the one hand,
\begin{align*}
{\mathbb P} & \left\{ Y_{\zeta_u} + \Pi_{\zeta_u} {\cal M}_{u} > u \right\}\\[.1cm]
 &={\mathbb P} \left\{ Y_{\zeta_u} + \Pi_{\zeta_u} {\cal M}_{u} > u , |Y_{\zeta_u}| \le \Delta(u) u\right\}
+{\mathbb P} \left\{ Y_{\zeta_u} + \Pi_{\zeta_u} {\cal M}_{u} > u, \: |Y_{\zeta _u}| > \Delta(u)u \right\}\\[.1cm]
&\leq  {\mathbb P} \left\{  \Pi_{\zeta_u} {\cal M}_{u} > \left(1 - \Delta(u) \right)u \right\}
+ {\mathbb P} \left\{ |Y_{\zeta_u}| > \Delta(u) u \right\};
\end{align*}
while on the other hand,
\begin{align*}
{\mathbb P} & \left\{ \Pi_{\zeta_u} {\cal M}_{u} > \left( 1 + \Delta(u) \right)u \right\}\\[.1cm]
&= {\mathbb P} \left\{ \Pi_{\zeta_u} {\cal M}_{u} > \left( 1 + \Delta(u) \right)u , |Y_{\zeta_u}| \le \Delta(u) u\right\} + {\mathbb P} \left\{ \Pi_{\zeta_u} {\cal M}_{u} > \left( 1 + \Delta(u) \right)u, \: |Y_{\zeta _u}| > \Delta(u)u \right\}\\[.1cm]
& \leq {\mathbb P} \left\{ Y_{\zeta_u} + \Pi_{\zeta_u} {\cal M}_{u} > u \right\}
+ {\mathbb P} \left\{ |Y_{\zeta_u}| > \Delta(u) u \right\}.
\end{align*}
Then, in view of \eqref{new-th2.1.3},
\begin{align}
{\mathbb P} & \left\{ \Pi_{\zeta_u} {\cal M}_{u} > \left( 1 + \Delta(u) \right)u \right\}-o \left( \frac{u^{-I(\tau)}}{\sqrt{\log u}} \right) \nonumber\\
& \leq {\mathbb P} \left\{ Y_{\zeta_u} + \Pi_{\zeta_u} {\cal M}_{u} > u \right\}
 \leq {\mathbb P} \left\{  \Pi_{\zeta_u} {\cal M}_{u} > \left(1 - \Delta(u) \right)u \right\}+o \left( \frac{u^{-I(\tau)}}{\sqrt{\log u}} \right).
\end{align}
Now apply \eqref{new-th2.1.2} to the left- and right-hand sides of this equation.  This yields
that
$${\mathbb P} \left\{ Y_{\zeta_u} + \Pi_{\zeta_u} {\cal M}_{u} > u \right\} \sim {\mathbb P} \left\{ \Pi_{\zeta_u} {\cal M}_{u}> u \right\} \ \mbox{ as } u \to \infty. $$
Hence the required result follows from \eqref{new-th2.1.2} and \eqref{pf-th2.1.1}.

\medskip

 \textbf{\textit{Step 2.}}
It remains to show that this constant $C(\tau)$  is positive and finite, and that the limit in this equation actually exists.

 \textbf{\textit{Step 2a.}}
First we prove existence of the limit. For this purpose
 we utilize the
$\alpha$-shifted measure defined previously in \eqref{measchange}.  Namely observe that by \eqref{eq: perpetuity} and \eqref{eq: max},
\[
 \frac{1}{\lambda^n(\alpha)} {\mathbb E}  \left[ M_n^\a \right]
 = {\mathbb E}_\alpha  \left[ \left( \max_{0 \le k \le n} Y_k \right)^\alpha \Pi_n^{-\alpha} \right]
= {\mathbb E}_\alpha \left[ \left( \max_{1 \le k \le n} \sum_{j=1}^k \widetilde{B}_j \left(\widetilde{A}_{j+1} \cdots \widetilde{A}_n\right)\vee 0  \right)^\alpha \right],
\]
where  $\widetilde{A}_j := 1/A_j$ and  $\widetilde{B}_j := B_j/A_j$ for all $j$.
By exchanging indices in this last expression, where we let $j \mapsto n+1-j$ in the expectation on the right-hand side, we then obtain
{
\begin{equation} \label{constP1}
\frac{1}{\lambda^n(\alpha)} {\mathbb E}  \left[ M_n^\a \right]
  = {\mathbb E}_{\alpha }\left [ \left( \max_{1\le k\le n} \sum_{j=k}^n \left( \widetilde{A}_1 \cdots \widetilde{A}_{j-1} \right) \widetilde{B}_j \vee 0 \right)^\alpha \right].
  \end{equation}
  Note that in this expression, the pair $(\widetilde{A},\widetilde{B})$ satisfies the following moment conditions:
\begin{align} \label{calculations-1}
&{\mathbb E}_{\alpha } \big[\log \widetilde{A} \big]  = -{\mathbb E}_{\alpha } \big[\log A \big] =- \frac 1{\lambda(\alpha)} {\mathbb E}  \big[A^{\alpha }\log A \big]<0; \nonumber\\[.2cm]
&{\mathbb E}_{\alpha} \big[ \widetilde A^{\alpha } \big]  = \frac{1}{\lambda (\alpha )}<1; \quad \mbox{and} \quad
{\mathbb E}_{\alpha} \big[ |\widetilde B |^{\alpha } \big]  =\frac{1}{\lambda (\alpha)} {\mathbb E} \big[ |B|^{\alpha} \big] <\infty.
\end{align}

To further analyze the limit in \eqref{constP1} as $n \to \infty$,
we first show:

\begin{Assert}\label{elemnt}
{\it Let $s_n = \sum_{j=1}^n d_j$ be an absolutely convergent series. Then the sequence
$$m_n = \max\big\{ d_n, d_{n-1}+d_n, \ldots ,d_1+ \cdots +d_n   \big\}$$ converges.
}
\end{Assert}

\noindent
{\bf Proof of the Assertion.}
It is sufficient to prove that $m_n$ is a Cauchy sequence. Fix $\epsilon>0$. Since the series is absolutely  convergent, there exists $N$ such that $\sum_{j>N}|d_j|<\epsilon $. Note
$$m_N = \max\big\{ d_N, d_{N-1}+d_N,\ldots,d_1+ \cdots +d_N   \big\},$$
and for any $p>N$,
$$m_p = \max\big\{ d_p, d_{p-1}+d_p, \ldots, d_{N+1}+ \cdots +d_p, \ldots,  d_1+ \cdots +d_N+d_{N+1}+ \cdots +d_p   \big\}.$$
Note that $m_p$ contains all of the factors that appear in $m_N$, but they  are modified by adding $d_{N+1}+ \cdots +d_p$ (which is at most $\epsilon $ in absolute value). Moreover, $m_p$ contains $N-p$ additional terms, but all of them are bounded, in absolute value, by $\epsilon $.
Therefore $\big|m_N - m_p\big|<\epsilon$ and $m_n$ is convergent.
\halmos

Now, in view of \eqref{calculations-1}, the perpetuity
$$
\widetilde Y_n = \sum_{j=1}^n \widetilde A_1\ldots \widetilde A_{j-1} \widetilde B_j
$$
converges $\P_\a$-a.s. Hence by the last assertion,
$$
X_n = \max_{1\le k\le n} \sum_{j=k}^n \left( \widetilde{A}_1 \cdots \widetilde{A}_{j-1} \right) \widetilde{B}_j \vee 0
$$ also converges $\P_\a$-a.s.  Set $X = \lim_{n \to \infty} X_n$.
Now $ X_n$ can be dominated by $R = \sum_{j=1}^\8 \widetilde A_1\ldots
\widetilde A_{j-1}|\widetilde B_j|$; and in view of \eqref{calculations-1} and \eqref{eq: moments}, we have that  $\E \left[ R^\a \right] <\8$.  Therefore, by the dominated convergence theorem,
$$
\lim_{n\to\8} \frac{1}{\lambda^n(\alpha)} {\mathbb E}  \left[ M_n^\alpha \right]
= \lim_{n\to\8} \E_{\alpha}\big[  X_n^{\alpha} \big] = \E_{\alpha}\big[  X^{\alpha} \big],
  $$ and this last expectation is finite.
This proves the existence of the limit.

 \textbf{\textit{Step 2b.}}
Finally, we prove that this limit is strictly positive.  To this end,
consider $\widetilde{Y}_n \vee 0$  as $n \to \infty$ (which we recognize { as a single term} in the maximum on the right-hand side of \eqref{constP1}).
Clearly, $\widetilde Y_n \le  X_n$.  Furthermore,
 $\widetilde{Y}_n$ converges to $\widetilde{Y}$ with $\E \big[  |\widetilde{Y} | ^{\alpha} \big] <\infty$, and
\begin{equation}
\lim _{n\to \infty} \frac 1{\lambda^n(\alpha)} \E\Big[ \left( \max\{ 0,Y_n \} \right)^\alpha \Big]= \lim_{n\to \infty }\E_{\alpha}\Big[ \left( \max\{ 0,\widetilde{Y}_n \} \right)^\alpha \Big]=\E_{\alpha}\Big[
  \left( \max\{ 0, \widetilde{Y} \} \right)^\alpha \Big].
\end{equation}
Also, observe that
$$ \E_{\alpha}\Big[
  \left( \max\{ 0, \widetilde{Y} \} \right)^\alpha \Big] \le \E_{\alpha}\left[
   X^\alpha \right].
$$

We claim that if ${\mathbb P} \left\{ A>1, \,B>0 \right\} >0$, then this last expectation is strictly positive.  Let $\widetilde{\pi}$ denote the probability law of $\widetilde{Y}$, and assume the assertion to be false.
Then $\widetilde Y \le 0$ $\P_\a$-a.s.; that is, $\supp (\widetilde{\pi} ) \subset (-\infty,0]$.
Notice that $\supp (\widetilde{\pi})$ must be $\nu _{\alpha }$-invariant a.s.\ under the action of $(\widetilde A, \widetilde B)$.   Also note that ${\mathbb P}  \left\{A>1, \,B>0 \right\}>0$ implies that
${\mathbb P}_{\alpha} \{ \widetilde{A} < 1,\, \widetilde{B} >0 \} >0$. Let $x_0 = \sup \big\{ x: x\in \supp (\widetilde{\pi}) \big\}$. Then $x_0 \leq 0$, but taking a pair $(\widetilde A, \widetilde B)$ such that $\widetilde A <1, \,\widetilde B >0$, we obtain
that $\widetilde A x_0 +\widetilde B > x_0$, and we are led to a contradiction.

This shows that the constant $C(\tau)$ in \eqref{pf-th2.1.10a} must be positive, thereby completing the proof of the theorem.
\halmos

\setcounter{equation}{0}
\section{Proof of Theorem \ref{thm:22}}
\label{section: 4}
\subsection{Preliminary considerations}
As in the previous section, define
$\overline{T}_u = (\log u)^{-1} \inf\{ n: \overline{Y}_n > u \},$
where $\overline{Y}_n = \sum_{i=1}^n \Pi_{i-1} |B_i|$.
First we establish an analog of Lemma \ref{lem:3.2} for the case $\tau = \rho$.
\begin{Lm}\label{lem:4.1}
 Assume that $\Lambda(\xi + \eta) <\8$ and $\Lambda_B(\xi + \eta) <\8$   for some $\eta>0$.
Then there exists a finite constant $\overline{D}$ and positive constant $\delta \equiv \delta(\eta)$ such that for all $u \ge 0$,
 \begin{equation} \label{COR-3.1.1}
 {\mathbb P} \left\{ \overline{T}_u \le \rho - L_\rho(u) \right\} \leq \overline{D} u^{-\xi}(\log u)^{-\delta }.
 \end{equation}
where $L_\rho(u) = b \sqrt{\{ \log (\log u)\}/\log u }$ for any constant $b >  \rho \left\{ 2 (\xi +1)+ \rho \sigma^2(\xi) \right\}$.
\end{Lm}

\Pf
Let
$\zeta_u = \lfloor \log u \left( \rho - L_\rho(u)\right) \rfloor$, then by definition
\begin{equation} \label{cor-3.1.1aNEW}
{\mathbb P} \left\{ \overline{T}_u \le \rho - L_\rho(u) \right\} =
{\mathbb P} \left\{ \overline{Y}_{\zeta_u} > u \right\}.
\end{equation}
Now apply Lemma \ref{lem:3.1} with $\alpha \equiv \xi$.  Since $\Lambda(\xi)=0$, it suffices to show that for some
$\epsilon \equiv \epsilon(u)$,
\begin{equation} \label{lm-4.1.1}
\overline{C}_{\zeta_u} u^{-\epsilon(u)} \leq \overline{D}(\log u)^{-\delta }.
\end{equation}
Let
\begin{equation} \label{lm-4.1.2}
\epsilon(u) =  \left( \frac{\log\,(\log u)}{\log u}  \right)^{1/2}.
\end{equation}

To analyze $\overline{C}_{\zeta_u}$, first note by \eqref{def-alpha} and the definition of $\rho$ that  $\Lambda^\prime(\xi) = \rho^{-1}$.  Hence for some finite constant
$D$,
\[
\overline{C}_{\zeta_u} \le D \left( \log u \right)^{2(\xi + \epsilon(u))+1} \exp\left\{ \log u( \rho - L_\rho(u)) \left(\frac{\epsilon(u)}{\rho} + {\epsilon^2(u)}\sigma^2(\xi) \right) \right\}.
\]
Thus, for sufficiently large $u$,
\begin{equation} \label{lm-4.1.3}
\overline{C}_{\zeta_u} u^{-\epsilon(u)}\leq D\exp \left\{ 2(\xi + 1) \log\,(\log u)  - \frac{1}{\rho} \log u  \Big( L_\rho(u) \epsilon(u) \Big) +\rho \log u \left( \epsilon^2(u) \sigma^2(\xi)\right) \right\}.
\end{equation}
Substituting the definitions of $L_\rho(u)$ and $\epsilon(u)$ into this last equation  yields
\begin{equation} \label{july31-11}
\overline{C}_{\zeta_u} u^{-\epsilon(u)}\leq D \exp \left\{ 2(\xi + 1) \log\,(\log u)  - \frac{b}{\rho} \log\,(\log u) + {\rho \sigma^2(\xi)} \log\,(\log u) \right\}
= D(\log u)^{-\delta },
\end{equation}
where $\delta > 0$ whenever $b >  \rho \left\{ 2 (\xi +1)+ \rho \sigma^2(\xi) \right\}$.  Thus we obtain \eqref{COR-3.1.1} for sufficiently large $u$ (with $\overline{D}= D$) and, hence,
with another choice $\overline{D} \ge D$, we obtain this equation for all $u \ge 0$.
\halmos

In the proofs below, it will be useful to observe that an analog of Lemma \ref{lem:4.1} also holds for the right tail of the hitting time of { $\{\overline  Y_n \}$ } to the level $u$.  To this end, set
\begin{align}
\label{eq: yn}
{ {\overline Y}^n } &= \sum_{k=n+1}^\infty \Pi_{k-1} |B_k|, \quad i =1,2,\ldots, \nonumber\\
{ {\overline T}^u }&=  (\log u)^{-1} \sup \left\{ n \in {\mathbb Z}_+:  { {\overline Y}^n } > u \right\}.
\end{align}

\begin{Lm}
 Assume that $\Lambda(\xi + \eta) <\8$ and $\Lambda_B(\xi + \eta) <\8$
 for some $\eta>0$.
 Then there are constants $C, \delta, b >0$ such that for every $u>e$,
 \begin{equation} \label{lm-4.2.1}
 {\mathbb P} \left\{ \overline{T}^u \ge \rho + L_\rho(u) \right\}\leq C u^{-\xi}(\log u )^{-\delta},
 \end{equation}
where $L_\rho(u) = b \sqrt{ \{ \log\,(\log u)\}/\log\, u }$.
\end{Lm}

\noindent
{\bf Proof.}
 Since $\sum_{k=1}^\infty k^{-2} = \pi^2/6$,  it follows that for some $\epsilon > 0$ (possibly dependent on $k$ and $u$),
\begin{equation}\label{lm-4.2.2}
{\mathbb P} \left\{\overline Y^{n}> u \right\} \le \sum_{k=n+1}^\infty {\mathbb P} \left\{ \Pi_{k-1} |B_{k}| > \frac{6u}{\pi^2 (k-n)^2} \right\}
\le \sum_{k=n+1}^\infty {\mathbb E} \Big[ \Pi_{k-1}^{\xi - \epsilon} |B_k|^{\xi-\epsilon} \Big] \left( \frac{\pi^2 (k-n)^2}{6u} \right)^{\xi - \epsilon}.
\end{equation}
Note by independence that
\[
{\mathbb E} \Big[ \Pi_{k-1}^{\xi - \epsilon} |B_k|^{\xi-\epsilon} \Big] = \left( {\mathbb E} \big[ A^{\xi - \epsilon} \big] \right)^{k-1} {\mathbb E} \big[ |B|^{\xi-\epsilon} \big] :=
   \left(  \lambda(\xi-\epsilon) \right)^{k-1} \lambda_B(\xi-\epsilon).
\]
Moreover, since $\Lambda(\xi)= 0,$ $\mu(\xi) = \rho^{-1}$, and $\Lambda$ is infinitely differentiable on the interior of its domain,
\[
\lambda(\xi- \epsilon) = e^{\Lambda(\xi- \epsilon)} \le \exp \left\{ - \frac{\epsilon}{\rho} + \frac{\epsilon^2 {\mathfrak l}}{2} \right\},
\]
where ${\mathfrak l} := \sup\{ \sigma^2(\alpha):  \xi - \epsilon \le \alpha \le \xi \}$.  Then using the continuity of $\sigma^2(\cdot)$, we have
that for sufficiently small $\epsilon$, ${\mathfrak l}/2 \le \sigma^2(\xi)$.
Hence, substituting the last two equations into \eqref{lm-4.2.2} yields
\begin{equation} \label{lm-4.2.3}
{\mathbb P} \left\{ {\overline Y^n} > u \right\} \le \left( \frac{\pi^2}{6} \right)^\xi u^{-\xi } \sum_{j=1}^\infty j^{2\xi} u^{\epsilon}\exp \left\{ (n+j-1) \left( - \frac{\epsilon}{\rho} + \epsilon^2 \sigma^2(\xi)
\right) \right\} \lambda_B(\xi - \epsilon).
\end{equation}

Now specialize to the case where $n \ge \log u(\rho + L_\rho(u))$.  Then with $\epsilon \equiv\epsilon (j) \equiv \epsilon(j,u)$, we obtain
\begin{align} \label{lm-4.2.4}
{\mathbb P} \left\{ \overline Y^n > u \right\} & \le \frac{\pi^2}{6} u^{-\xi} \sum_{j=1}^\infty j^{2\xi} \exp \Big\{-\frac{\epsilon(j)L_\rho(u)\log u}{\rho}
   -(j-1) \frac{\epsilon  (j)}{\rho} \nonumber\\[-.1cm]
& \hspace*{5.25cm}    + (n+j-1) \epsilon^2 (j) \sigma^2(\xi) \Big\}\lambda_B(\xi - \epsilon (j)).
\end{align}
Now choose
$$\epsilon(j) = \gamma \ \frac {  L_\rho(u)\log u+(j-1)} { \rho {\sigma}^2(\xi) (n+j-1)},$$
where $\gamma$ is a positive constant. Since this expression remains bounded as $u\to\infty$ (uniformly in $j\ge 1$), the constant
 $\gamma$ can be chosen such that $\epsilon(j)$ is arbitrarily small. Then for  $n(u) = \lfloor \log u(\rho + L_\rho(u)) \rfloor$,  $b\geq \rho$, and $\gamma_1 = \gamma - \gamma^2$, we obtain by \eqref{lm-4.2.4} that \begin{eqnarray} \label{lm-4.2.5}
{\mathbb P} \left\{ \overline Y^{n(u)} > u \right\} & \le & C u^{-\xi} \sum_{j=1}^\infty j^{2\xi} \exp \left\{ -\frac{\gamma_1 (L_\rho(u)\log u +j-1)^2}{4\rho^2 \sigma^2(\xi) (n(u)+j-1)} \right\} \nonumber\\[.2cm]
& \le & C  u^{-\xi}  \Bigg(\big( n(u) \big)^{2\xi +1}\exp \left\{ -  \gamma_1 b\log (\log u )/16\rho ^2{\sigma^2(\xi)}\right \} \nonumber\\[-.1cm]
 && \hspace*{1.75cm} + \sum_{j  \ge n(u)+1} j^{2\xi} \exp \left\{ -\gamma_1 (j-1)/8 \rho^2 {\sigma^2(\xi)} \right\}
  \Bigg) \nonumber\\[.2cm]
  &\leq & C u^{-\xi}(\log u)^{2\xi +1-\gamma_1 b\slash 16\rho ^2 {\sigma^2(\xi)}},
\end{eqnarray}
since for $j\leq n(u)$ we have
$$
\frac{(L_\rho(u)\log u +j-1)^2}{4\rho^2 {\sigma} (n(u)+j-1) }\geq \frac{b^2\log (\log u)}{8\rho ^2(\rho +b){\sigma}}\geq \frac{b\log (\log u)}{16\rho ^2{\sigma}}.$$
Thus \eqref{lm-4.2.1} follows from \eqref{lm-4.2.5} upon choosing $b \ge \max \{\rho, 16{\sigma^2(\xi)}  \rho ^2 ( 2 \xi +1)/\gamma_1\}$.
\halmos

From the previous lemma, we  draw two conclusions.  First, we observe that this lemma combined with Lemma \ref{lem:4.1} may be used to prove a strengthening of Lemma \ref{Lm2.1},
thus establishing a conditional law of large numbers for the scaled first passage time of $\{ Y_n \}$ to level $u$.

\begin{Lm}\label{lem:4.3}
Let $L_\rho(u)$ be given as in Lemma \ref{lem:4.1}, and assume that
$\Lambda$ and $\Lambda_B$ are finite in a neighborhood of $\xi$ and the law of $\log A$ is nonarithmetic.  Then
\begin{equation} \label{lm-4.3.1}
\lim_{u \to \infty} {\mathbb P} \left\{ \left. \left| T_u - \rho \right| \ge L_\rho(u) \: \right| \: T_u < \infty \right\} = 0.
\end{equation}
\end{Lm}

\Pf  Note $\overline{Y}_n \ge Y_n$, for all $n$, implying that $\{ \overline{T}_u \le \rho - L_\rho(u) \} \supset \{ T_u \le \rho - L_\rho(u) \}$.
Consequently, it follows by Lemma \ref{lem:4.1} that
 \begin{equation} \label{lm-4.3.2}
 {\mathbb P} \left\{ T_u \le \rho - L_\rho(u) | T_u < \infty \right\} = o (1) \quad \text{as} \quad u \to \infty.
 \end{equation}

Next, set $n_u = \lceil \log u (\rho + L_\rho(u) \rceil$ and define ${ R}_{n} = { M} - {M}_{n}$.
Observe that { ${ R}_n \le  {\overline Y}^{n} $, } for all $n$.  Hence by Lemma \ref{lem:4.1},
\begin{equation} \label{lm-4.3.3}
{\mathbb P} \left\{{  { R}_{n_u} } > u \right\} = o(u^{-\xi}) \quad \mbox{as} \:\: u \to \infty.
\end{equation}
Thus, by repeating the argument following
\eqref{new-th2.1.2} above,
we obtain that the tail decay of { ${ M} = {M}_{n_u} + {R}_{n_u}$ }
is dominated by the larger of the tails of { ${M}_{n_u}$ and ${R}_{n_u}$}, respectively, which must necessarily be the tail of { ${ M}_{n_u}$}; that is,
\begin{equation} \label{lm-4.3.4}
\lim_{u \to \infty} u^\xi {\mathbb P} \left\{ { M} > u \right\} = \lim_{u \to \infty} u^\xi {\mathbb P} \left\{ { M}_{n_u} > u \right\}.
\end{equation}
Since  $\{ { M}_{n_u} > u \} \subset \{ { M} > u \}$, it follows that
\begin{equation} \label{lm-4.3.5}
{\mathbb P} \left\{ T_u \ge \rho + L_\rho(u) | T_u < \infty \right\} =  1 - \frac{ {\mathbb P} \left\{ { M}_{n_u} > u \right\} }{ {\mathbb P} \left\{ { M} > u \right\} } =   o (1) \quad \text{as} \quad u \to \infty,
\end{equation}
as required.
\halmos

From the perspective of our main theorems, a more important consequence to be drawn from Lemma \ref{lem:4.1} is the convergence of a certain measure $H$ to Lebesgue measure.
We will establish this convergence in the Assertion given in the proof of Theorem \ref{thm:22}.

\subsection{Establishing the main result for the critical case}

{\bf Proof of Theorem \ref{thm:22}.}
 \textbf{\textit{Step 1.}}
Let $L_\rho(u) \ge b \sqrt{ \{\log\,(\log u)\}/\log u }$, where $b >  \rho \left\{ 2(\xi + 1) + \rho \sigma^2(\xi)/2 \right\}$.  Then by Lemma \ref{lem:4.1},
\begin{equation} \label{pf-thm2.2.1}
{\mathbb P} \left\{ \overline{T}_u \le \rho - L_\rho(u) \right\} = o ( u^{-\xi} ) \quad \text{as} \quad u \to \infty.
\end{equation}
Set
\[
\zeta_u = \lfloor \log u \left( \rho - L_\rho(u)\right) \rfloor, \quad \quad \rho_u =  \lfloor \rho \log u \rfloor,
\]
and
\begin{equation}
\label{eq: Mu}
{\cal M}_{u} = \max_{\zeta_u  < n \le \rho_u} \Big\{ B_{\zeta_u+1} + A_{\zeta_u+1} B_{\zeta_u+2} + \cdots
  + \left( A_{\zeta_u+1} \cdots A_{n-1} \right) B_n \Big\} \vee 0.
\end{equation}
Then arguing as in the proof of Theorem \ref{thm:21} (specifically, by repeating the argument following \eqref{new-th2.1.2}), we obtain by \eqref{pf-thm2.2.1} that
\begin{equation} \label{pf-thm2.2.2}
{\mathbb P} \left\{ T_u \le \rho \right\} = {\mathbb P} \left\{ \Pi_{\zeta_u} {\cal M}_{u} > u \right\} \left(1+ o(1)\right) \quad \mbox{as} \:\: u \to \infty.
\end{equation}

To analyze the right-hand side of this equation, we begin by observing, as in the proof of Theorem \ref{thm:21}, that
\begin{equation} \label{pf-th2.2.3}
{\mathbb P} \left\{ \Pi_{\zeta_u} {\cal M}_u > u \right\} = \int_{\reals } {\mathbb P} \left\{ \log \Pi_{\zeta_u} > \log u - s\right\} dF_u(s),
\end{equation}
where $F_u$ denotes the probability distribution function of $\log {\cal M}_u$.  Then apply Petrov's theorem to handle the probability on the
right-hand side.

First observe (cf.\ \eqref{new-P}, \eqref{pf-th2.1.4b}) that
\[
\log \Pi_{\zeta_u} > \log u - s \Longleftrightarrow \frac{S_{\zeta_u}}{\zeta_u} >  \frac{\log\,u - s}{\zeta_u} := \frac{1}{\rho} + \gamma_u,
\]
where, for a deterministic function $\bfdelta_u$ with $|\bfdelta_u| \le 1\slash \rho $, we have
\begin{equation} \label{pf-th2.2.4}
\gamma_u = \frac{1}{\zeta_u} \left( \frac{L_\rho(u)}{\rho} \log u - s + \bfdelta_u \right)\quad \mbox{and}\quad
\zeta_u \gamma_u^2 = \frac{1}{\zeta_u} \left( \frac{L_\rho(u)}{\rho} \log u - s + \bfdelta_u \right)^2.
\end{equation}
Now let $\Delta > 0$ and consider ${\mathbb P} \left\{ \Pi_{\zeta_u} M_u > u, {\cal G}_u \right\}$, where
\begin{eqnarray*}
{\cal G}_u &=& \left\{\omega \in \Omega:  \log \, M_u(\omega) \in  \left[0, D(u) + \Delta \sqrt{\zeta _u}\sigma (\xi ) \right]  \right\},\\[.1cm]
D(u) &=& \frac{L_\rho(u)}{\rho} \log u + \bfdelta _u.
\end{eqnarray*}

Note that when $s \in {\cal H}_u := \left[0,D(u) +\Delta \sqrt{\zeta _u}\sigma (\xi ) \right]$ (corresponding to the event ${\cal G}_u$ in \eqref{pf-th2.2.3}), we have
by elementary calculations that $\gamma_u \to 0$ and $\zeta_u \gamma_u^3 \to 0$ as $u \to \infty$, uniformly for $s \in {\cal H}_u$.  However, we do not have that
$\zeta_u \gamma_u^2 \to 0$ as $u \to \infty$.   Thus, focusing on the exponential term in Petrov's theorem, we see that the first- {\it and} second-order terms must be retained in
the expansion (in contrast to the proof of Theorem \ref{thm:21}, where it was sufficient to analyze the first-order term), while the third-order term may again be neglected.
Consequently, by Petrov's Theorem \ref{thm:3.1Petrov}, we obtain that
\begin{equation} \label{pf-th2.2.5}
{\mathbb P} \left\{ \Pi_{\zeta_u} {\cal M}_u > u, {\cal G}_u \right\} = \frac{1}{\xi\sigma(\xi) \sqrt{2\pi \zeta_u}} \int_0^{D(u) +\Delta \sqrt{\zeta _u}\sigma (\xi )}  { g(u,s)} dF_u(s) \left(
  1 + o(1) \right),
\end{equation}
where
\begin{equation} \label{pf-th2.2.6}
{ g(u,s)} = u^{-\xi}e^{\xi s}  \exp \left\{ - \frac{1}{2 \sigma^2(\xi) \zeta_u} \left( D(u) - s \right)^2  \right\}.
\end{equation}

Next introduce the transformation
\begin{equation} \label{trans-T}
{\mathbb T}_u(s) = \frac{1}{\sigma(\xi) \sqrt{\zeta_u}} \left( s - D(u) \right),
\end{equation}
and let $G_u(E) = F_u ({\mathbb T}_u^{-1}(E)),$ for all $E \in {\cal B}(\reals)$.
Then after a change of variables (\shortciteN{PB86}, p.\ 219), we obtain that
\begin{eqnarray} \label{pf-th2.2.7}
{\mathbb P} \left\{ \Pi_{\zeta_u} {\cal M}_u > u, {\cal G}_u \right\} & = & \frac{u^{-\xi}}{\xi\sigma(\xi) \sqrt{2\pi \zeta_u}} \int_{-D(u)/ \sigma(\xi) \sqrt{\zeta_u} }^{\Delta} e^{\xi \,{\mathbb T}^{-1}_u(z)}
   e^{-z^2/2}dG_u(z) \left(1 + o(1) \right) \nonumber\\[.2cm]
& = & \frac{u^{-\xi}}{\sqrt{2\pi}} \int_{-D(u)/ \sigma(\xi) \sqrt{\zeta_u} }^{\Delta } e^{-z^2/2} dH_u(z) \left( 1 + o(1) \right)
\end{eqnarray}
as $u \to \infty$, where for any set $E \in {\cal B}(\reals)$,
\begin{equation} \label{pf-th2.2.8a}
H_u (E) = \frac{1}{\xi \sigma(\xi) \sqrt{\zeta_u}} \int_E e^{\xi \,{\mathbb T}^{-1}_u(z)} dG_u(z) = \frac{1}{\xi \sigma(\xi) \sqrt{\zeta_u}} \int_{{\mathbb T}^{-1}_u(E)} e^{\xi s} dF_u(s).
\end{equation}

 \textbf{\textit{Step 2.}}
Our strategy now is to characterize the measure $H_u$ on $(-\infty , \Delta]$, and then to handle the remaining part (namely, ${\mathbb P} \left\{ \Pi_{\zeta_u} {\cal M}_u > u, {\cal G}_u^c \right\}$) by a separate argument.
This remaining part will turn out to be asymptotically negligible.
To characterize the measure $H_u$ on $(-\infty,\Delta]$, we establish the following.
\begin{Assert}
{\rm (i)}  Let ${\cal C}_M$ denote the constant appearing in \eqref{eq: asymptotic},
let $l$ denote Lebesgue measure on $\reals$, and let ${\bf C}_K$ denote the collection of continuous functions on $(-\infty,0)$ with compact support.  Then for any $f \in {\bf C}_K$,
\begin{equation} \label{pf-th2.2.8}
\lim_{u \to \infty} \int_{(-\infty,0)} f(s) dH_u(s) = {\cal C}_M \int_{(-\infty,0)}f(s) dl(s).
\end{equation}\\
\hspace*{.1in} {\rm (ii)}  There is a constant $\overline{\cal C}$ such that for every $-\infty <v<w< \infty$ and every $u  \ge 0$,
\begin{equation} \label{pf-th2.2.9}
H_u(v,w) \le \overline{\cal C} (w-v)+ \frac{K}{\sqrt{\zeta _u}}, \quad \mbox{for some } K < \infty.
\end{equation}
\end{Assert}

\noindent
{\bf Proof of the Assertion.}
We begin with the proof of (ii).
  Let $-\infty < v < w < \infty$, and set $v^\ast(u) = {\mathbb T}_u^{-1}(v),$ $w^\ast (u) = {\mathbb T}_u^{-1}(w)$, and  let $\overline{F}_u(z) = 1- F_u(-\infty,z)$, $z \in \reals$.
Then from \eqref{pf-th2.2.8a} and an integration by parts,
\begin{equation} \label{pf-th2.2.10}
H_u(v,w) = - \frac{1}{\xi \sigma(\xi) \sqrt{\zeta_u}} \left\{ e^{\xi w^\ast(u)}\overline{F}_u \left(w^\ast(u) \right) - e^{\xi v^\ast(u)}\overline{F}_u \left(v^\ast(u) \right) \right\} +
  \frac{1}{\sigma(\xi) \sqrt{\zeta_u}} \int_{v^\ast(u)}^{w^\ast(u)} e^{\xi s} \overline{F}_u(s) ds.
\end{equation}

To analyze the first term on the right-hand side of this equation, observe that the distribution of ${\cal M}_u$ is stochastically dominated by $M$.
Moreover, as we have observed in Section \ref{section: 2},
\begin{equation} \label{pf-th2.2.12}
\lim_{u \to \infty} u^\xi {\mathbb P} \left\{ M > u \right\} = {\cal C}_M.
\end{equation}
Consequently, there exists a finite positive constant $\overline{\cal C}$ such that
\begin{equation} \label{pf-th2.2.13}
{\mathbb P} \left\{ { M} > u \right\} \le \overline{\cal C} u^{-\xi}, \quad \mbox{for all} \:\: u \ge 0.
\end{equation}
From \eqref{pf-th2.2.13}, we obtain that $e^{\xi s} \overline{F}_u (s) \le \overline{\cal C}$ for all $s \in \reals$.  Hence we conclude that the first term on the right-hand side of \eqref{pf-th2.2.10} is dominated by $O\left(1/\sqrt{\zeta _u} \right)$ as $u \to \infty$, where $\zeta_u \to \infty$, independent of the choice of $v$ and $w$.

For the second term on the right-hand side of this equation, first note that by
 inverting the function ${\mathbb T}_u$ defined in \eqref{trans-T}, we obtain
\begin{equation} \label{inverseT}
{\mathbb T}^{-1}_u (t) = D(u) + t \sigma(\xi) \sqrt{\zeta_u}, \quad t \in \reals.
\end{equation}
Hence
\[
w^\ast(u) - v^\ast(u) := {\mathbb T}^{-1}_u(w) - {\mathbb T}^{-1}_u(v) = \sigma(\xi) \sqrt{\zeta_u}(w-v),
\]
and consequently
\begin{equation} \label{july31-1}
  \frac{1}{\sigma(\xi) \sqrt{\zeta_u}} \int_{v^\ast(u)}^{w^\ast(u)} e^{\xi s} \overline{F}_u(s) ds\leq \frac{\overline{\cal C}}{\sigma(\xi) \sqrt{\zeta_u}}\left(w^\ast(u)-v^\ast(u) \right)=\overline {\cal C}(w-v).
\end{equation}
Thus we have established (ii).

Turning to the proof of (i), now suppose that $-\infty < {\mathfrak a} \le v < w \le {\mathfrak b} < 0$.
We begin by observing that
\begin{equation} \label{pf-th2.2.14}
\lim_{u \to \infty} e^{\xi s} \overline{F}_u(s) = {\cal C}_M, \quad \mbox{uniformly for} \:\: s \in {\mathbb T}^{-1}_u([{\mathfrak a},{\mathfrak b}]).
\end{equation}

To establish this claim, first recall that $F_u$ is the distribution function of $\log {\cal M}_u$, where ${\cal M}_u$ was defined in \eqref{eq: Mu}.
Set
\begin{equation} \label{def-n.u}
m_u = \rho_u - \zeta_u := \lfloor \rho \log u \rfloor - \lfloor \log u (\rho - L_\rho(u) ) \rfloor = L_\rho(u)\log u + \bfdelta^\ast_u,\quad |\bfdelta^\ast_u| \le 1.
\end{equation}
Now to prove \eqref{pf-th2.2.14}, we use Lemma \ref{lem:4.1}.
Namely, we show  that
for  some finite constant $\overline{D}$ and some positive constant $\delta$,
\begin{equation}\label{pf-th2.2.14c}
e^{\xi s}{\mathbb P} \left\{ { {\overline Y}^{m_u} }  > e^s \right\} \leq \overline{D} s^{-\delta }, \quad s \in {\mathbb T}^{-1}_u\left([{\mathfrak a},{\mathfrak b}]\right),
\end{equation}
 where $\overline Y^n$ was defined in \eqref{eq: yn}.

Before establishing \eqref{pf-th2.2.14c}, we first observe that \eqref{pf-th2.2.14c} implies \eqref{pf-th2.2.14}.
For this purpose,
let $\{ { M}_n\}$ and $\{ { R}_n \}$ be defined as in the proof of Lemma \ref{lem:4.3}, and observe that these definitions imply that ${ R}_{m_u} \le
{ {\overline Y}^{m_u} }$ and
${ M} = { M}_{m_u} + {R}_{m_u}$. Therefore, by \eqref{pf-th2.2.14c},
\begin{equation}\label{pf-th2.2.14d}
e^{\xi s}{\mathbb P} \left\{{ R}_{m_u} > e^s \right\} \leq \overline{D} s^{-\delta }.
\end{equation}
Arguing as in the proof of Lemma \ref{lem:4.3}, we then conclude that
\begin{equation}
\lim _{s\to \infty }e^{\xi s} {\mathbb P} \left\{ { M}_{m_u} > e^s \right\} = \lim _{s\to \infty }e^{\xi s} {\mathbb P} \left\{ { M} > e^s \right\} = {\cal C}_M.
 \end{equation}
Moreover, as a straightforward consequence of the definitions, we see that $\log M_{m_u}$ is equal in distribution to $\log {\cal M}_u$, which has the distribution $F_u$.
Hence \eqref{pf-th2.2.14} follows.

To establish \eqref{pf-th2.2.14c}, in view of Lemma \ref{lem:4.1},  it is enough to observe that
$$
m_u>(\rho +L_{\rho }(e^s))s.$$
[We apply Lemma \ref{lem:4.1} to $e^s$ instead of $u$.]
Now
\begin{align} \label{july30}
&\rho s <  \rho {\mathbb T}^{-1}({ b}) =  L_{\rho }(u)\log u +\rho \bfdelta _u +\rho { b} \sigma (\xi )\sqrt{\zeta _u},\nonumber\\[.1cm]
& \hspace*{1cm} m_u = L_{\rho }(u)\log u +\bfdelta^\ast_u,
\end{align}
where $|\rho \bfdelta |, |\bfdelta^\ast|\leq 1$. Hence we need to show that for sufficiently large $u$,
\begin{equation}\label{pf-th2.2.14e}
L_{\rho }(e^s)s+2< -\rho { b} \sigma (\xi )\sqrt{\zeta _u}, \quad \mbox{where} \:\: {\mathfrak b} < 0.
\end{equation}
But for sufficiently large $u$, we see from the first equation in \eqref{july30} that (as $L_\rho(u) \ge b \sqrt{ \{\log\,(\log u)\}/\log u }$) we have
$$s\leq \frac{b}{\rho}\sqrt{\log (\log u)}\sqrt{\log u}.$$
Hence
$$
L_{\rho }(e^s)s\leq b\sqrt{s\log s}=o(\sqrt{\zeta _u}) \quad \mbox{as}\:\: u \to \infty,$$ and \eqref{pf-th2.2.14e} follows.
Thus we have established \eqref{pf-th2.2.14c} and consequently \eqref{pf-th2.2.14}.

Now returning to \eqref{pf-th2.2.10} and focusing on the second term on the right-hand side of this equation,
observe by the uniform convergence in \eqref{pf-th2.2.14} that
\begin{equation} \label{pf-th2.2.15a}
\lim_{u \to \infty} \frac{1}{\sigma(\xi) \sqrt{\zeta_u}} \int_{v^\ast(u)}^{w^\ast(u)} e^{\xi s} \overline{F}_u(s) ds = \lim_{u \to \infty} \frac{{\cal C}_M(w^\ast(u) - v^\ast(u))}{\sigma(\xi) \sqrt{\zeta_u}} ={\cal C}_M(w-v),
\end{equation}
where the last step follows as in \eqref{july31-1}.  Since the first term on the right of \eqref{pf-th2.2.10} is $O(1/\sqrt{\zeta_u})$, as we have shown in the proof of (i), we conclude that
\begin{equation} \label{july31-2}
\lim_{u \to \infty} H(u,v) = {\cal C}_M (w-v), \quad - \infty < {\mathfrak a} \le v < w \le {\mathfrak b} < 0.
\end{equation}
Then taking Riemann sums, we obtain that for any $f \in {\bf C}_K$,
\begin{equation} \label{pf-th2.2.17}
\lim_{u \to \infty} \int_{-\infty}^0 f(s) dH_u(s) = {\cal C}_M \int_{-\infty}^0 dl(s),
\end{equation}
as required.
\halmos

\textbf{\textit{Step 3.}}
Now returning to the proof of the main theorem, we split our interval into three parts,
\[
\big[-\frac{D(u)}{\sqrt{\zeta _u}\sigma (\xi)},-J\big], \quad \big[ -J,-\Delta \big],\quad \big[-\Delta , \Delta \big], \quad \mbox{ where } 0 < J < \infty,
\]
and observe by part (i) of the previous assertion that for any $J$,
\[
\lim _{u\to \infty} \int _{-J}^{-\Delta }e^{-z^2/2}dH_u(s) = {\cal C}_M \int_{-J}^{-\Delta }e^{-z^2/2} dl(s),
\]
while by part (ii) of the assertion,
\begin{eqnarray*}
 \int _{-D(u)/\left(\sigma (\xi)\sqrt{\zeta _u}\right)}^{-J}e^{-z^2/2}dH_u(s) & \leq & \overline {\cal C} \int_{-\infty}^{-J}e^{-z^2/2} dl(s)\leq \overline {\cal C} e^{-J^2/2},\\[.2cm]
 \int _{-\Delta }^{\Delta} e^{-z^2/2}dH_u(s) & \leq & \overline {\cal C} \int_{-\Delta }^{\Delta }e^{-z^2/2} dl(s)\leq 2\overline {\cal C}\Delta .
\end{eqnarray*}
Letting $J \to \infty$ yields
\begin{equation} \label{pf-th2.2.18}
\lim_{u \to \infty} u^\xi {\mathbb P} \left\{ \Pi_{\zeta_u} {\cal M}_u > u, {\cal G}_u \right\}
= \frac{{\cal C}_M}{\sqrt{2\pi}} \int_{-\infty}^{\Delta} e^{-z^2/2} dl(z) = \frac{{\cal C}_M}{2} + o(1) \quad \mbox{as} \:\: \Delta \to 0.
\end{equation}

 \textbf{\textit{Step 4.}}
It remains to show that the restriction to ${\cal G}_u$ can be removed on the left-hand side of this last equation.

 \textbf{\textit{Step 4a.}}
We begin by removing the restriction that $\log {\cal M}_u \le D(u) +\Delta \sqrt{\zeta _u}\sigma (\xi )$.   To this end, letting ${\mathbb T}_u$ be defined as in \eqref{trans-T}, then it is sufficient to show that
for any $\Delta > 0$,
\begin{equation} \label{pf-th2.2.20}
{\mathbb P}  \left\{ \Pi_{\zeta_u} {\cal M}_u > u, \: \log {\cal M}_u > {\mathbb T}_u^{-1}(\Delta ) \right\} = o(1) u^{-\xi} \quad \mbox{as} \:\: u \to \infty.
\end{equation}

Let $\Delta>0$ be given, and set ${\cal G}_{u,k} = \left\{\omega \in \Omega:  \log {\cal M}_u(\omega) - {\mathbb T}_u^{-1}(\Delta ) \in (k-1,k]\right\}$, $k=1,2,\ldots.$
Then
\[
\left\{  \Pi_{\zeta_u} {\cal M}_u > u, \: \log {\cal M}_u > {\mathbb T}_u^{-1}(\Delta ) \right\} = \bigcup_{k \in \pintegers} \left\{  \Pi_{\zeta_u} {\cal M}_u > u, \: {\cal G}_{u,k} \right\}.
\]
Moreover,
\begin{equation} \label{pf-th2.2.21}
{\mathbb P} \left\{  \Pi_{\zeta_u} {\cal M}_u > u, \: {\cal G}_{u,k} \right\} \le {\mathbb P} \left\{ \Pi_{\zeta_u} >  u e^{-\left({\mathbb T}_u^{-1}(\Delta ) + k \right)} \right\} {\mathbb P} \left\{ {\cal G}_{u,k} \right\}
  \le  u^{-\xi} e^{\xi {\mathbb T}_u^{-1}(\Delta )} e^{\xi k} {\mathbb P} \left\{ {\cal G}_{u,k} \right\}
\end{equation}
by Chebyshev's inequality.  In addition, by applying Lemma \ref{lem:3.1} with $\epsilon \equiv \epsilon(u) = (\log u)^{-1/3}$, we obtain that
\begin{equation} \label{pf-th2.2.22}
{\mathbb P} \left\{ {\cal G}_{u,k} \right\}= {\mathbb P} \left\{ {\cal M}_u > e^{{\mathbb T}_u^{-1}(\Delta )} e^{k-1} \right\} \le \overline{C}(u) e^{-(\xi + \epsilon(u)) {\mathbb T}_u^{-1}(\Delta )} e^{-(\xi +\epsilon(u))(k-1)},
\end{equation}
where $\overline{C}(u)$ corresponds to the quantity $\overline{C}_n$ appearing in the statement of Lemma \ref{lem:3.1}.  To identify the growth rate of this function as $u \to \infty$,
note that  ${\cal M}_u \stackrel{\cal D}{=} \left\{ Y_i:  0 \le i \le \rho_u - \zeta_u \right\}$, where (from the definitions following \eqref{pf-thm2.2.1}) we have that
$m_u :=\rho_u - \zeta_u = L_\rho(u) \log u + \bfdelta^\ast_u$ for $|\bfdelta^\ast_u| \le 1$; cf.\ \eqref{def-n.u}.  Thus, in Lemma \ref{lem:3.1}, we must replace the parameter $n$ with
$m_u$.  Since $\mu(\xi) { =} \rho^{-1}$, we consequently obtain that, for some positive constant $K$,
\[
\overline{C}(u) \le K \left( L_\rho(u) \log u \right)^{2(\xi + 1)} \exp \left\{ \epsilon(u)  (L_\rho(u)/\rho) \log u \right\}.
\]
(The term ``$(n-1) \epsilon^2 \sigma^2(\a)$'' of Lemma \ref{lem:3.1} is negligible here, since $m_u \sim L_\rho(u) \log u$, which grows at rate $\sqrt{\log\, (\log u)} \sqrt{\log u}$,
while $\epsilon(u)$ is chosen so that it decays at rate $(\log u)^{-1/3}$.)
Moreover (cf.\ \eqref{inverseT}),
\[
{\mathbb T}_u^{-1}(\Delta ) = \frac{L_\rho(u) \log u}{\rho} + \Delta \sigma(\xi) \sqrt{\zeta_u}+\bfdelta _u.
\]
Combining these last two equations yields
\begin{equation} \label{pf-th2.2.23}
\overline{C}(u) e^{-\epsilon(u){\mathbb T}_u^{-1}(\Delta )} \le K^\prime \left(L_\rho(u) \log u\right)^{2 (\xi + 1)} e^{-\epsilon(u)\Delta \sigma(\xi) \sqrt{\zeta_u}} = O\left(\exp \left\{ -\Delta \sigma(\xi)(\log u)^{1/6} \right\}\right)
\end{equation}
as $u \to \infty$.  Hence by \eqref{pf-th2.2.22},
\begin{equation} \label{july31-4}
{\mathbb P} \left\{ {\cal G}_{u,k} \right\} = O\left(\exp \left\{ -\Delta \sigma(\xi)(\log u)^{1/6} \right\}\right) e^{-\xi {\mathbb T}_u^{-1}(\Delta)}  e^{-(\xi +\epsilon(u))(k-1)} \quad \mbox{as} \:\: u \to \infty.
\end{equation}
Substituting this equation into \eqref{pf-th2.2.21}, we conclude that {
\begin{align} \label{pf-th2.2.24}
{\mathbb P} & \left\{ \Pi_{\zeta_u} {\cal M}_u > u, \: \log {\cal M}_u > {\mathbb T}_u^{-1}(\Delta ) \right\} = \sum_{k=1}^\8 {\mathbb P} \left\{  \Pi_{\zeta_u} {\cal M}_u > u, \: {\cal G}_{u,k} \right\} \\
 &\qquad= O \left(\exp \left\{ -\Delta \sigma(\xi)(\log u)^{1/6} \right\} \right) u^{-\xi} \sum_{k=0}^\infty e^{-k \epsilon(u)} \nonumber\\[.2cm]
 & \qquad = O\left(\exp \left\{ -\Delta \sigma(\xi)(\log u)^{1/6}\right\} \right) u^{-\xi}(\log u)^{1/3}= o(1) u^{-\xi} \quad \mbox{as} \:\: u \to \infty,
\end{align} }
which establishes \eqref{pf-th2.2.20}.

 \textbf{\textit{Step 4b.}}
Finally, observe by Chebyshev's inequality followed by a Taylor expansion that
\begin{align} \label{pf-th2.2.25}
{\mathbb P} \left\{ \Pi_{\zeta_u} {\cal M}_u > u, \: \log {\cal M}_u < 0 \right\} \le {\mathbb P} \left\{ \Pi_{\zeta_u}  > u \right\}
  \le u^{-\xi - \epsilon(u)} \exp \left\{ \zeta_u \left( \frac{\epsilon(u)}{\rho} + { \epsilon^2(u) \sigma^2(\xi)} \right) \right\}.
\end{align}
Next recall that  $\zeta_u := \lfloor \log u \left( \rho - L_\rho(u) \right) \rfloor$; thus,
\[
u^{-\epsilon(u)} e^{ \epsilon(u) \zeta_u/\rho} \le \exp \left\{-\epsilon(u) \left( \frac{L_\rho(u)}{\rho} \right)\log u + 1 \right\}.
\]
Now choose $\epsilon(u) = \left(\log u \right)^{-1/2}$.  Then on the right-hand side of the previous equation, the exponential term tends to $-\infty$ as $u \to \infty$.  Moreover, with
this choice of $\epsilon(u)$, we also have that
$\zeta_u \epsilon^2(u)$ is bounded as $u \to \infty$.  Thus we conclude that
\begin{equation} \label{pf-th2.2.26}
{\mathbb P} \left\{ \Pi_{\zeta_u} {\cal M}_u > u, \: \log {\cal M}_u < 0 \right\} = o(1) u^{-\xi} \quad \mbox{as} \:\: u \to \infty,
\end{equation}
as required.

 \textbf{\textit{Step 5.}}  It remains to prove that the constant ${\cal C}_M$ is strictly positive.
 To this end, let $Y$ be defined as in Section 2.1, let $\pi$ denote the probability law of $Y$, and let ${\cal C}_Y$ be given as in \eqref{eq: asymptotic}.
 Clearly ${\cal C}_Y \le {\cal C}_M$.
 Thus, it is sufficient to analyze the case where ${\cal C}_Y=0$. Then, by a result of \shortciteN{GL2014},
   the support of $\pi$ is unbounded
 from below. 

 We claim that Collamore and Vidyashankar (2013b, Remark 2.3 and Section 9) can be applied
to obtain the representation formula \eqref{JCAVconst1}; hence ${\cal C}_M$ is strictly positive.
For this purpose, we need to show that Lemma 5.1 (iii) of that paper holds without the continuity assumption given
there (namely, condition ($H_0$) of that paper).
Let the forward process $\{ M_n^\ast \}$ be defined as in \eqref{eq: recursions}, and  let $P$ denote its transition kernel.
Then $P$ satisfies a minorization condition, namely,
\begin{equation}
P(x,E) \ge {\bf 1}_{\mathbb C}(x) \eta(E), \qquad x \in \reals, \: E \in {\cal B}(\reals),
\end{equation}
where ${\mathbb C} = \{ 0 \}$ and $\eta$ denotes the probability law of $B^+$.
Note that this minorization condition is nontrivial, since $\supp(\pi) \cap (-\infty,0] \not= \emptyset$ implies that $Y_n^\ast = A_n Y_{n-1}^\ast + B_n$ hits
$(-\infty,0]$ with positive probability; and hence, $\{ M_n^\ast \}$ hits $\{0 \}$ with positive probability, as these two processes agree up until the first time that
$Y_n^\ast \le 0$.
Thus, in particular, $\{ M_n^\ast \}$ is $\psi$-irreducible (\shortciteN{EN84}, Remark 2.1).

Moreover, if $B^+$ has a density on some subinterval of $(0,\infty)$, then the set $[0,K]$ is petite.   Indeed, letting $\pi^+$ denote the stationary measure of $\{M_n^\ast \}$,
then $\supp(\pi^+) \supset \supp(\eta)$ (the support of $B^+$), implying that $\supp(\pi^+)$ is of second category.
Hence, since $[0,K]$ is a compact set and
$\{ M_n^\ast \}$ is a weak Feller chain, we may apply Remark 2.7 (i) of
\shortciteN{ENPT82} to conclude that $[0,K]$ is petite.

But if $B^+$ does not have a density on some subinterval of $(0,\infty)$, then we can dominate $\{M_n^\ast\}$ from below by the process $\{\widetilde{M}_n^\ast\}$, where
\[
\widetilde{M}_n^\ast = \left(A_n \widetilde{M}_{n-1}^\ast + \widetilde{B}_n \right)^+,  \qquad \widetilde{B}_n = B_n + \zeta_n,
\]
where $\{ \zeta_n \}$ is an i.i.d.\ sequence, independent of $\{ (A_n,B_n) \}$, such that $\zeta$ has a smooth density supported on the interval $(-\delta,0)$ for some $\delta > 0$.
The process $\{\widetilde{M}_n^\ast\}$ is regularly varying at infinity with parameter $\xi$; that is, it satisfies \eqref{eq: asymptotic} and
 the above argument can be applied  to conclude that the corresponding  constant ${\cal C}_{\widetilde{M}}$ is positive.  Since $\{M_n^\ast\}$ dominates the process $\{\widetilde{M}_n^\ast\}$, we  conclude that ${\cal C}_M$ is also strictly positive.
\halmos

\setcounter{equation}{0}
\section{Proof of Theorems \ref{thml1} and \ref{thml2}}
\label{section: 5}
\subsection{Proof of Theorem \ref{thml1}}
Before we proceed with the proof of Theorem \ref{thml1}, it is worthwhile to observe that there exists a measure satisfying the hypotheses of this theorem and, in particular,
 \eqref{eq: count1}.

\begin{lem}  There exists a measure satisfying the assumptions of Theorem \ref{thml1}.
\end{lem}

\noindent
{\bf Proof.}
Take an arbitrary measure $\nu$ with continuous nonvanishing density $h$ and such that, if $A$ has law $\nu$, then $\E[ \log A]<0$ and $\E [A^{\a+\epsilon}] <1$ for some $\a >1$.
(For example, one could choose $h(a)=C\eta ^{-1}(1+Ca)^{-\a -2\eps-1}$ with $C$ sufficiently large and $\eta = \int _0^{\8}(1+u)^{-\a -2\eps-1}\ du $.)  Define the family of probability measures
$$
\nu_t(da) = t h(ta) da, \quad t > 0.
$$
Let $\Lambda_t(\b) = \log \E_{\nu_t} \left[ A^\b \right]$ denote the cumulant generating function.
Then
\[
\Lambda_t(\a) = \log\bigg( \int_{\reals} a^\a h(ta) tda  \bigg) = -\a \log t + \Lambda(\a),
\]
and
\[
\mu_t(0) := \Lambda'_t(0) =  \int_\reals \log a\; h(ta) tda = -\log t + \mu(0).\\
\]
Hence
$$
\mu_t(0) - \Lambda_t(\a) = (\a-1)\log t + \mu(0) - \Lambda(\a),
$$
and choosing $t$ appropriately large, we have that $\mu_t(0) > \Lambda_t(\a)$. Thus the measure $\nu_t$ satisfies the hypotheses of Theorem \ref{thml1}.
\halmos

The proof of Theorem \ref{thml1} will now be based on the following two lemmas.  In these lemmas, we study the joint distribution of $\Pi_n := A_1 \cdots A_n$
and $Y_{n} := 1 + A_1 + \cdots + (A_1 \cdots A_{n-1})$ as $n \to \infty.$

%
%
%
%

\begin{lem}
\label{lem: 5.9}
Let $\beta \in {\rm int} (\dom \Lambda)$ be chosen such that $\Lambda(\b) < \infty$, and set  $\tau_\b = (\mu(\b))^{-1}$ and $k_u = \lfloor \log u/\mu(\b)\rfloor$.
Then there are positive constants $D_0, \g, a,b$
with $\g b<1$,
such that for sufficiently $\mbox{large }u$,
\begin{equation} \label{aaa1}
\P\big\{ \g u\le \Pi_{k_u-1} \le \g a u, \quad \! Y_{k_u} \le \g b u \big\} \ge \frac{D_0}{\sqrt{\log u}} u^{-I(\tau_\beta)}.
\end{equation}
\end{lem}

\begin{lem}
\label{lem: 5.99}
Assume $\E \left[ |\log A|^3 \right] <\8$, and set $\epsilon(m) = e^{(m-1)\mu(0)}$. Then there are positive constants $D_1,c,d$ such that for sufficiently large $m$,
\begin{equation} \label{aaa2}
\P \left\{ \epsilon(m) \le \Pi_{m-1} \le c \eps(m), \quad \! Y_m \le d  \right\}  \ge \frac {D_1}{\sqrt m}.
\end{equation}
\end{lem}

Heuristically, Lemma 5.2 and 5.3 can be understood as follows.  Since $\mu(\beta) = {\mathbb E}_\beta \left[ \log A \right]$, it follows that in the $\beta$-shifted measure, the random
walk
\[
S_n = \log \Pi_n := \sum_{i=1}^n \log A_i, \quad n=1,2,\ldots
\]
will reach the boundary at level $\log u$ at approximately the time $\tau_\beta \log u$, i.e., roughly at time $k_u$.  Hence it follows from standard large deviation arguments (based on the Bahadur-Rao approximation, cf.\
\shortciteN{ADOZ93}) that
\begin{equation} \label{abc1}
\P\left\{ \log u + \log \gamma \le S_{k_u}\le  \log u + \log (\g a) \right\} \sim \frac{D_0(\tau)}{\sqrt{\log u}} u^{-I(\tau_\beta)}
\end{equation}
as $u \to \infty$.   (In particular, \eqref{abc1} can be concluded from Petrov's Theorem 3.1, stated above.)
Hence, \eqref{aaa1} states that on $\left\{ \log u + \log \gamma \le S_{k_u}\le  \log u +  \log(\g a) \right\}$, we have with
high probability that
$\log Y_{k_u} \le  \log u + \log(\gamma b)$ for some positive constant $b$.  The latter event can be expected, since
$Y_n := 1 + A_1 + \cdots + (A_1 \cdots A_{n-1})$,
and, in the $\beta$-shifted measure, this process will grow as $u \to \infty$ since $\mu(\beta) > 0$.  Roughly speaking, $\{ Y_n \}$ will then be dominated by its last term,
namely $\Pi_{n-1}$, where $\log \Pi_{n-1} = S_{n-1}$.

In a similar way, Lemma 5.3 can be viewed, roughly speaking, as a consequence of the Berry-Essen theorem, which studies the asymptotic behavior of $\{S_n\}$ around its central tendency, namely
around $n\mu(0) = n {\mathbb E} \left[ \log A \right]$.  Then the  Berry-Essen theorem yields the estimate \eqref{aaa2}, but without the ``$Y_m \le d$'' term on the left-hand side.
Note that when $\{ S_n \}$ follows a trajectory which is close it its mean trajectory, one expects $\Pi_n \downarrow 0$ and $n \to \infty$ and hence $\{ Y_n \}$ to be convergent.

Nonetheless, the formal proofs of Lemmas 5.2 and 5.3 are quite technical and hence postponed to the end of this section.
Before turning to these rigorous proofs, we first show how our main result may be deduced from these two lemmas.

\medskip

%
%
%
%

\noindent
{\bf Proof of Theorem \ref{thml1}.}  Let $\alpha \equiv \alpha(\tau)$ (where $\alpha(\tau)$ is given as in \eqref{def-alpha}), and let $0 < \alpha < \beta < \xi$ be chosen such that
$$
\Lambda'(\b) = \frac{\Lambda'(\a)}{p}.
$$
Later, we will choose $\b$ close to $\alpha$ and $p$ close to one, but the precise choices of these constants will be fixed only at the end of the proof.
Note that since $\mu(\cdot) := \Lambda^\prime(\cdot)$,
\[
\tau = \frac{1}{\mu(\alpha)} = \frac{1}{p \mu(\beta)} := \frac{\tau_\beta}{p}.
\]
For $q = 1-p$,  now define
\[
n_u = \left\lfloor \tau \log u \right\rfloor, \quad
\quad k_u = \lfloor pn_u \rfloor, \quad m_u  = n_u-k_u.
\]

We begin by writing
\begin{equation}
Y_{n_u} = Y_{k_u} + \Pi_{k_u-1}A_{k_u}Y^\prime_{m_u},
\end{equation}
where
$$
Y'_{m_u} = 1+A^\prime_1+ \cdots + A^\prime_1 \cdots A^\prime_{m_u-1} \quad \mbox{for} \quad A^\prime_i = A_{n_u+i}.
$$
Note that $Y^\prime_{m_u}$
is independent of $Y_{k_u}, \Pi _{k_u-1}$ and $A_{k_u}$. Let
\[
\Pi'_{m_u-1} = A^\prime_1 \cdots A^\prime_{m_u-1},\quad
 \bfepsilon_u = e^{(m_u-1) \mu(0)},
\]
and
$$
\Omega_u = \left\{ \g u\le \Pi_{k_u-1} \le \g a u, \quad \! Y_{k_u} \le b \g u, \quad \!
\bfepsilon_u \le \Pi_{m_u-1}^\prime \le c \bfepsilon_u, \quad \! Y_{m_u}^\prime \le d \right\},
$$
where the constants $a,b,c,d$ are given as in Lemmas  \ref{lem: 5.9} and \ref{lem: 5.99}.
Applying Lemmas \ref{lem: 5.9} and \ref{lem: 5.99}, we then conclude that there exists a constant $D$ such that, for sufficiently large $u$,
\begin{equation}
\P \left\{ \Omega_u \right\} \ge \frac{D}{\sqrt{m_u}\sqrt{\log u}}  u^{-I(\tau_\b)}.
\end{equation}

Next observe
\begin{align}
&\P\left\{ Y_{n_u-1}\le u \mbox{ and } Y_{n_u}>u \right\}\nonumber\\[.2cm]
&\hspace*{1cm} = \P\left\{ Y_{k_u} + \Pi_{k_u-1}A_{k_u} Y_{m_u-1}^\prime  \le u \mbox{ and }
Y_{k_u} + \Pi_{k_u-1}A_{k_u}Y^\prime_{m_u} > u \right\} \nonumber\\[.2cm]
&\hspace*{1cm} \ge \P\left\{  \left( \frac{u-Y_{k_u}}{\Pi_{k_u-1}Y'_{m_u}} < A_{k_u} \le
 \frac{u-Y_{k_u}}{\Pi_{k_u-1}Y'_{m_u-1}} \right)   \cap\ \Omega_u \right\} \nonumber\\[.2cm]
&\hspace{1cm} =  \P \left\{  \left( C_- \le  \frac{u-Y_{k_u}}{\Pi_{k_u-1}Y'_{m_u}} < A_{k_u} \le
 \frac{u-Y_{k_u}}{\Pi_{k_u-1}Y'_{m_u-1}}\ \le  C_+ \right) \cap\ \Omega_u \right\}
\end{align}
for certain constants $C_+$ and $C_-$, since on the set $\Omega_u$ we have
\[
C_+ := \frac{1}{\gamma} \ge \frac{u}{\Pi_{k_u-1} Y_{m_u-1}^\prime}
 \qquad\mbox{and}\qquad C_- :=  \frac{1- \gamma b}{\gamma a d} \le \frac{u - Y_{k_u-1}}{\Pi_{k_u-1} Y_{m_u}^\prime} .
\]

Notice that on the set $\Omega_u$, we have that for sufficiently large $u$,
$$
 \frac{u-Y_{k_u}}{\Pi_{k_u-1}Y^\prime_{m_u-1}}
- \frac{u-Y_{k_u}}{\Pi_{k_u-1}Y^\prime_{m_u}}  \ge \frac{(u-Y_{k_u})\Pi^\prime_{m_u-1}}{\Pi_{k_u-1} d^2/2} \ge \frac{1-\gamma b}{\gamma a d^2/2}\cdot \bfepsilon_u := \gamma^\ast \bfepsilon_u.
$$
Therefore, for sufficiently large $u$, $A_{k_u}$ must belong to a random interval of length at least $\gamma^\ast \bfepsilon_u$.

\medskip

Let $I_{\eps_u}$ be an arbitrary interval of length $\g^\ast \bfepsilon_u$. Since the density of $A$ is bounded from below on the interval $[C_+,C_-]$ by some constant $\d$, we have
\begin{eqnarray} \label{aab0}
\P\left\{ Y_{n_u-1}\le u \mbox{ and } Y_{n_u} >u \right\} &\ge& \inf_{I_{\bfepsilon_u}\subset [C_-,C_+]} \P\left\{ \big( A_{k_u}\in I_{\bfepsilon_u} \big) \cap \Omega_u \right\}\nonumber\\[.1cm]
&\ge&  \inf_{I_{\bfepsilon_u}\subset [C_-,C_+]} \P\left\{  A \in I_{\bfepsilon_u} \right\} \P\left\{ \Omega_u \right\}\nonumber\\[.1cm]
&\ge&  \frac{\d \bfepsilon_u}{\sqrt{m_u}} \frac{D}{\sqrt{\log u}}  u^{-I(\tau_\b)}.
\end{eqnarray}
We emphasize that in this computation, $\Omega_u$ has been defined so that $A_{k_u}$ is independent of $\Omega_u$.
We now elaborate on the last term.  Our objective is to compare the decay rate $\bfepsilon_u u^{-I(\tau_\beta)}$ to the ``expected" decay rate governed by $u^{-I(\tau)}.$
To this end, note by an application of \eqref{eq: rate2} (cf.\ \eqref{cor-3.1.2a}) that
\begin{equation}\label{exponent}
u^{ I(\tau) - I(\tau_\b)} \bfepsilon_u =\exp\Big\{\a \log u-\Lambda(\a)n_u^\ast -\b\log u+\Lambda(\b) k_u^\ast+m_u^\ast \mu(0) \Big\},
\end{equation}
where $n_u^\ast = \tau \log u$, $k_u^\ast = p n_u^\ast$ $(= \tau_\beta \log u)$ and $m_u^\ast = n_u^\ast - k_u^\ast$.
Now estimate the exponent in \eqref{exponent} from below.  Recalling that $\mu(\alpha) = \tau^{-1}$,
we obtain
\begin{align} \label{aab1}
\a \log u &-\Lambda(\a)n_u^\ast -\b\log u+\Lambda(\b) k_u^\ast + m_u^\ast \mu(0)\nonumber\\[.1cm]
&=n_u^\ast\Big\{ \mu(\a)( \a - \b)  + p(\Lambda(\b) - \Lambda(\a)) + q(\mu(0) -\Lambda(\a))    \Big\} \nonumber\\[.1cm]
&=p\Big( \Lambda(\b)-\Lambda(\a) -\mu(\a)(\b-\a) \Big) +q\mu(\a)(\a-\b) + q(\mu(0) - \Lambda(\a) \big) \nonumber\\[.1cm]
&\ge q \Big(\mu(\a) (\a -\b )+\mu(0)-\Lambda (\a ) \Big),
\end{align}
since $\Lambda(\b)-\Lambda(\a) -\mu(\a)(\b-\a) = \Lambda^{\prime\prime}(\theta) > 0$, for some $\theta \in [\a,\b]$.

Since we are assuming that $\mu(0) > \Lambda(\a)$, we see that when $\b$ is close to $\a$,
 the last expression in \eqref{aab1} is strictly positive.  Thus, $\bfepsilon_u u^{-I(\tau_\beta)}$ decays at a slower polynomial rate than $u^{-I(\tau)}$.
 Hence the required result follows from \eqref{aab0}.
 \halmos

We now return to the proofs of Lemmas 5.2 and 5.3.
%


\medskip

\noindent
{\bf Proof of Lemma 5.2.}
 \textbf{\textit{Step 1.}}
By Theorem 3.1, there exists a constant $C_0$ such that
\begin{equation}
\label{eq: wt1}
\P\big\{\gamma u \le \Pi_{k_u-1} \le \g a u  \big\} =\frac{C_0}{\sqrt{\log u}} \gamma^{-\beta} u^{-{I(\tau_\b)}} (1+o(1)) \quad \mbox{as} \quad u \to \infty.
\end{equation}
The main step is to prove
 that there exist positive constants $\epsilon$ and $C_1$ and a constant $0<\d<1$ such that, for sufficiently large $u$,
\begin{equation}
\label{eq: t2}
\P\left\{\g u \le \Pi_{k_u-1}, \quad\! \Pi_{k_u-i-1} > \frac{\g b u}{2i^2}  \right\} \le  \frac{C_1}{\sqrt{\log u}}  \left( \g^{-\b} b^{-\epsilon} \d^i \right) u^{-I(\tau_\b)}
 \quad \mbox{ for } i=1,\ldots, k_u-1.
\end{equation}
The final result follows easily from \eqref{eq: wt1} and \eqref{eq: t2}, since then we obtain
\begin{align}
\P \Big\{ \g u & \le \Pi_{k_u-1} \le \g a u, \quad \! Y_{k_u} \le \g b u \Big\} \nonumber\\[.1cm]
&\ge \P\bigg\{ \g u \le \Pi_{k_u-1} \le \g a u, \quad \! \Pi_{k_u-i-1} \le \frac{\g b u}{2 i^2}\ \mbox{for all } i \in \{1,\ldots, k_u-1 \}\bigg\}\nonumber\\[.1cm]
&\ge \P\Big\{ \g u \le \Pi_{k_u-1} \le \g a u\Big\} - \sum_i \P \bigg\{ \g u \le \Pi_{k_u-1} \le \g a u, \quad \! \Pi_{k_u-i-1} > \frac{\g b u}{2 i^2}\bigg\} \nonumber\\[.1cm]
&\ge \bigg( C_0(1+o(1))- \frac { b^{-\epsilon} C_1}{1 - \delta} \bigg)\frac{\g^{-\b}}{\sqrt{\log u}}  u^{-I(\tau_\b)} \quad \mbox{as} \quad u \to \infty.
\end{align}
The result then follows after choosing the constant $b$ sufficiently large.

\medskip

 \textbf{\textit{Step 2.}}  We now prove \eqref{eq: t2}. For this purpose, we consider two cases, namely the case $i \le K \log k_u$ and then the case
 $i > K \log k_u$, where $K$ is a large positive constant and $k_u = \lfloor \tau_\beta \log u \rfloor.$

 \textbf{\textit{Case 1:}}  First assume that $i \le K \log k_u$, and suppose that a constant $L$ has been chosen such that
  \begin{equation} \label{abc10}
  -L \beta -\Lambda(\beta)K + \frac{1}{2} \le - \eta < 0.
  \end{equation}
Clearly,
\begin{align} \label{abd1}
\P\bigg\{\g u &\le \Pi_{k_u-1}, \quad \! \Pi_{k_u-i-1}\ge \frac{\g b u}{i^2}  \bigg\} \nonumber\\[.1cm]
& \le \P\bigg\{ \Pi_{k_u-i-1}\ge { \g b u \cdot k_u^L}  \bigg\}
+ \P\bigg\{ \g u \le \Pi_{k_u-1}, \quad \! \frac{ \g b u }{i^2}\le \Pi_{k_u-i-1}\le {\g b  u \cdot k_u^L}  \bigg\}.
\end{align}
Now the first term on the right-hand side is asymptotically negligible, since it follows by an application of Chebyshev's inequality that for some finite constant $C_2$,
\begin{eqnarray}
{\mathbb P} \left\{ \Pi_{k_u-i-1} \ge \gamma b u \cdot k_u^L \right\} &\le& \left( \g b u  \right)^{-\beta} k_u^{-L \beta} \left(\lambda(\beta) \right)^{k_u - i-1} \nonumber\\[.2cm]
& \le & \frac{C_2}{\sqrt{\log u}} u^{-I(\tau_\beta)} \left( k_u^{-L \beta + \frac{1}{2}} e^{-i \Lambda(\beta)} \right) \gamma^{-\b},
\end{eqnarray}
and by \eqref{abc10} and $\Lambda(\beta) < 0$, we have that for all $i \le K \log k_u$,
\[
k_u^{-L \beta + \frac{1}{2}} e^{-i \Lambda(\beta)} \le e^{-\eta \log k_u } \searrow 0 \quad \mbox{as} \quad u \to \infty.
 \]

 Thus, it is sufficient to focus on the second term on the right-hand side of \eqref{abd1}.  To this end, first note that
 \begin{align}
 \P\bigg\{ \g u & \le \Pi_{k_u-1}, \quad \! \frac{\g b u }{i^2}\le \Pi_{k_u-i-1}\le {\g b u \cdot k_u^L}  \bigg\}\nonumber\\[.1cm]
& \le \sum_{0 \le l \le \log(i^2 k_u^L)} \P\bigg\{ \frac{\g b u }{i^2}\cdot e^l \le \Pi_{k_u-i-1}< \frac{\g b u }{i^2}\cdot e^{l+1}
 \bigg\} \P\bigg\{ \Pi_i \ge \frac{i^2}{b e^{l+1}} \bigg\} \nonumber\\[.1cm]
& \le \sum_{0 \le l \le (L+1) \log k_u} \P\bigg\{  \Pi_{k_u-i-1}\ge \frac{\g b u }{i^2}\cdot e^{l}
 \bigg\} \P\bigg\{ \Pi_i \ge \frac{i^2}{b e^{l+1}} \bigg\}
\end{align}
for sufficiently large $u$.
The strategy is then to estimate the first term on the right-hand side by Petrov's theorem, and to estimate the second term using Chebyshev's inequality.
Using that $S_n := \log \Pi_n$, we see that the first term can be written as
\[
\P \left\{ S_{k_u-i-1}\ge \log \left( \frac{\g b u}{i^2} \right) + l \right\},
\]
where $k_u := \lfloor \tau \log u \rfloor \Longrightarrow \left\{ \log \left( \g b u/i^2\right)  + l \right\}/(k_u-i-1)  \sim \tau^{-1}$ as $u \to \infty$.
Note that since $i \le K \log k_u$, the conditions of Theorem 3.1 are easily verified.
Then by an application of Petrov's Theorem 3.1, we obtain that
\begin{equation}
\P \left\{ S_{k_u-i-1} \ge \log \left( \frac{\g b u}{i^2} \right) + l \right\} = \frac{1}{\beta  \sigma(\beta) \sqrt{2 \pi k_u}}
 \left( \frac{\gamma b e^l}{i^2} \right)^{-\beta} u^{-\beta} e^{(k_u-i)\Lambda(\b)}
\left(1 + o(1) \right)
\end{equation}
as $u \to \infty$, uniformly for $i \le K \log u$.  Recalling that $u^{-\beta}
 \exp \left\{\tau_\beta \log u  \cdot \Lambda(\beta)\right\} = u^{-I(\tau_\beta)}$ (cf.\ \eqref{cor-3.1.2a}), we then obtain that for some finite constant $C_3$,
\begin{equation}
 \P\left\{  \Pi_{k_u-i-1}\ge \frac{\g b u }{i^2}\cdot e^{l} \right\} \le \frac{C_3}{\sqrt{\log u}} u^{-I(\tau_\beta)} \cdot i^{2\beta} \gamma^{-\beta} b^{-\beta} e^{-l \beta} e^{-i \Lambda(\beta)},
 \qquad u \ge \mbox{some }U_0.
\end{equation}
Moreover, by Chebyshev's inequality, we have that for $\epsilon > 0$ sufficiently small
\begin{equation}
\P\left\{ \Pi_i \ge \frac{i^2}{b e^{l+1}} \right\} \le i^{-2(\b-\eps)} b^{\beta - \epsilon} e^{(l+1)(\beta-\epsilon)} e^{i \Lambda(\beta - \epsilon)},
\end{equation}
since $i^{2(\beta - \epsilon)} \ge 1$.  Next observe that $\Lambda(\beta) - \Lambda(\beta - \epsilon) = \epsilon \Lambda^\prime(\bar{\beta}) >0 $ for some $\bar{\beta} \in (\beta - \epsilon,
\beta)$,  where positivity of $\Lambda^\prime$ follows from the convexity of $\Lambda$.
   Hence we obtain from the previous two equations that, for some positive constant $C_4$ and sufficiently large $u$,
\begin{align} \label{abe1}
 \sum_{0 \le l \le (L+1) \log k_u} &\P\bigg\{  \Pi_{k_u-i-1}\ge \frac{\g b u }{i^2}\cdot e^{l}
 \bigg\} \P\bigg\{ \Pi_i \ge \frac{i^2}{b e^{l+1}} \bigg\}  \nonumber\\[.1cm]
 &\le \frac{C_4}{\sqrt{\log u}} \left(  \gamma^{-\beta} b^{-\epsilon} \right) u^{-I(\tau_\beta)}  i^{2\eps} e^{-i \epsilon \Lambda^\prime(\bar{\beta})}
   \sum_{l=0}^\infty e^{-\epsilon l},
\end{align}
which yields \eqref{eq: t2}.

\medskip

\textbf{\textit{Case 2:}}   Now suppose that $i > K \log k_u$.  Then by Chebyshev's inequality,
\begin{align}
\P\bigg\{\g u & \le \Pi_{k_u-1}, \quad\! \Pi_{k_u-i-1} > \frac{\g b u}{2i^2}  \bigg\} \nonumber\\[.2cm]
& \le \sum_{l=0}^\infty \P\bigg\{ \frac{\g b u }{i^2}\cdot e^l \le \Pi_{k_u-i-1}< \frac{\g b u }{i^2}\cdot e^{l+1}
 \bigg\} \P\bigg\{ \Pi_i \ge \frac{i^2}{b e^{l+1}} \bigg\} \nonumber\\[.1cm]
& \le \sum_{l=0}^\infty \left( \frac{ \gamma b u}{i^2}  \cdot e^l \right)^{-\beta} e^{(k_u - i-1) \Lambda(\beta)} \cdot  \left( \frac{i^2}{b e^{l+1}} \right)^{-(\beta- \epsilon)} e^{i \Lambda(\beta-\epsilon)}
  \nonumber\\[.1cm]
& \le \frac{1}{\lambda(\beta)}  \left( \gamma^{-\beta} b^{-\epsilon} \right) u^{-I(\tau_\beta)} \left( i^{2\epsilon} e^{-i \epsilon \Lambda^\prime(\bar{\beta})} \right) \sum_{l=0}^\infty e^{-l \epsilon},
\end{align}
where $\bar{\beta}$ is given as in \eqref{abe1}.  Hence, using that $(\log u)^{-1/2} \sim \sqrt{\tau} \exp\{\log k_u/2\} \le \sqrt{\tau} e^{i/K}$ for $i > K \log u$, we then
obtain that for some positive constant $C_5$,
\begin{align}
\P\bigg\{\g u & \le \Pi_{k_u-1}, \quad\! \Pi_{k_u-i-1} > \frac{\g b u}{2i^2}  \bigg\} \nonumber\\[.2cm]
& \le \frac{C_5}{\sqrt{\log u}} \left( \gamma^{-\beta} b^{-\epsilon} \right) u^{-I(\tau_\beta)} \left( i^{2\epsilon} e^{-i \epsilon \Lambda^\prime(\bar{\beta})} e^{i/K} \right),
\end{align}
which establishes \eqref{eq: t2} upon choosing $K$ sufficiently large.
\halmos

\noindent
{\bf Proof of Lemma 5.3.}    \textbf{\textit{Step 1.}}   With a slight abuse of notation, we will write $\epsilon_m$ in place of $\epsilon(m)$ throughout the proof.

First we prove that there exist finite constants $c$ and ${\mathfrak M}$ such that
\begin{equation} \label{abf1}
\P\left\{ \eps_m \le \Pi_{m-1} \le c \eps_m \right\} \ge \frac{\log c}{2 \sigma \sqrt{2\pi}} \frac{1}{\sqrt{m}}, \qquad \mbox{for all } m \ge {\mathfrak M},
\end{equation}
where $\s^2 = {\rm Var}(\log A)$.

Since $\mu(0) = {\mathbb E} [X_1]$, it follows by the Berry-Esseen theorem (\shortciteN{VP95}, Theorem 5.5) that for all $m$,
\begin{equation}
\label{eq: berry-essen}
\sup_x \bigg| \P\bigg\{ \frac{S_m - m \mu(0)}{\s \sqrt m} < x\bigg\} - \Phi(x) \bigg| \le \frac{{\cal A} \E \left[ |X_1-\E X_1|^3\right]}{\s^3}\frac 1{\sqrt m}:=\frac{\rho }{\sqrt m},
\end{equation}
where $S_m := \sum_{j=1}^m \log A_i$ and $\Phi$ denotes the normal distribution function, and where ${\cal A}$ is a universal constant.
Hence for any $c > 1$,
\begin{equation}
\P\Big\{ 0 \le S_m - m \mu(0)  \le \log c \Big\} \ge  \left( \Phi \left( \frac{\log c}{\s \sqrt{m}} \right) - \Phi(0) \right) - \frac{2 \rho}{\sqrt{m}}.
\end{equation}
Now it follows from the definitions that $\log \left\{ \Pi_{m-1}/\epsilon_m \right\} = S_{m-1} - (m-1) \mu(0)$.   Thus, from the previous equation we obtain that
\begin{equation}
\P\left\{ \eps_m \le \Pi_{m-1} \le c \eps_m \right\}
\ge \frac{\log c}{\s \sqrt{2\pi}} \frac{1}{\sqrt{m-1}} e^{-\frac{1}{2} \left(\frac{\log c}{\s \sqrt{m-1}}\right)^2} - \frac{2\rho}{\sqrt{m-1}}.
\end{equation}
Then choosing $c$ sufficiently large yields \eqref{abf1}.

\smallskip

\textbf{\textit{Step 2.}}  Next we show that for sufficiently large $m$,
\begin{equation} \label{abf2}
\P\left\{ \eps_m \le \Pi_{m-1} \le c \eps_m, \quad \! Y_m > d \right\} \le \frac{{\cal B}}{d^\theta \sqrt{m}},
\end{equation}
where ${\cal B}$ and $\theta$ are finite positive constants.  Noting that
\[
\P\left\{ \eps_m \le \Pi_{m-1} \le c \eps_m, \quad \! Y_m \le d \right\} \le \sum_j \P\left\{ \eps_m \le \Pi_{m-1} \le c \eps_m, \quad \! \Pi_j > \frac{d}{2j^2} \right\},
\]
we then divide the sum on the right-hand side into two parts.

Let $p < 1$, and suppose that $\theta > 0$ has been chosen such that $\Lambda(\theta) < 0$.  Then, on the one hand,
\begin{equation} \label{abf3}
\sum_{j>pm} \P\left\{ \Pi_j > \frac{d}{2 j^2} \right\} \le
\sum_{j>pm} \frac{(2j^2)^\theta}{d^\theta}  e^{j \Lambda(\theta)} = o \left( \frac{1}{\sqrt{m}} \right)
\end{equation}
as $m \to \infty$.  On the other hand, we also have
\begin{align} \label{abf4}
\sum_{j \le pm} \P\bigg\{ \eps_m & \le \Pi_{m-1} \le c \eps_m, \quad \! \Pi_j > \frac{d}{2j^2} \bigg\} \nonumber\\[.1cm]
& \le \sum_{j \le pm} \sum_{k \ge 0} \P \left\{\frac{d}{2j^2} \cdot e^k \le \Pi_j < \frac{d}{2j^2} \cdot e^{k+1} \quad \mbox{and} \quad \eps_m  \le \Pi_{m-1} \le c \eps_m
  \right\} \nonumber\\[.1cm]
& \le \sum_{j \le pm} \sum_{k \ge 0} \P \left\{ \Pi_j \ge \frac{d}{2j^2} \cdot e^k \right\} \P \left\{ \frac{2 \eps_m j^2}{d e^{k+1}} < \Pi_{m-j} \le \frac{2 c \eps_m j^2}{d e^k}
   \right\}.
\end{align}
To estimate the last quantity on the right-hand side, apply once again the Berry-Esseen theorem, noting that $j \le pm \Longrightarrow m-j > (1-p)m$.  This yields (after a short
computation) that
\begin{equation}
\P \left\{ \frac{2 \eps_m j^2}{d e^{k+1}} < \Pi_{m-j} \le \frac{2 c \eps_m j^2}{d e^k} \right\}
 \le \frac{1}{\sqrt{2\pi}} \int_{\frac{C_{j,k}}{\s\sqrt{m-j}}}^{\frac{C_{j,k}+\log c+1}{\s\sqrt{m-j}}} e^{-\frac{x^2}2}dx + \frac {2\rho}{\sqrt m} \le \frac{{\cal B}^\prime}{\sqrt m},
\end{equation}
where $C_{j,k}=j \mu(0) + \log(2j^2) - \log d - k - 1$ and ${\cal B}^\prime$ is a finite constant.  Note that this integral is taken over an interval of length $(\log c +1)/(\sigma \sqrt{m-j})$.
Consequently, with $\theta$ chosen as before and ${\cal B}^{\prime\prime}$ a positive constant, we obtain that for sufficiently large $m$,
\begin{equation} \label{abf20}
\sum_{j \le pm} \P\bigg\{ \eps_m  \le \Pi_{m-1} \le c \eps_m, \quad \! \Pi_j > \frac{d}{2j^2} \bigg\}  \le  \sum_{j \le pm} \sum_{k \ge 0}  \frac{(2 j^2)^\theta}{d^\theta e^{\theta k}} e^{j \Lambda(\theta)}
  \frac{{\cal B}^{\prime\prime}}{\sqrt{m}}  \le \frac{{\cal B}}{d^\theta \sqrt{m}}
\end{equation}
for some positive constant ${\cal B}$, as required.

\smallskip

\textbf{\textit{Step 3.}}   Finally observe that if $d$ is chosen sufficiently large in the previous equation, then the decay in \eqref{abf1} dominates the decay
in \eqref{abf20}.  Consequently, the required result follows from \eqref{abf1} and \eqref{abf2}.
\halmos

\medskip

Finally, we remark that
Theorem 2.3 also holds with $B > 0$ a.s.\ (but not necessarily constant, as was assumed in the previous proofs).  However, in this case, the proofs become noticeably more technical.
Thus, in order to emphasize the main ideas in the proofs, we have restricted our attention to the case $B=1$.

To prove Theorem 2.3 for $B > 0$ a.s., we would need Lemma 5.2 at the required level of generality,
and also Lemma 5.3 slightly modified.  Namely, in place of Lemma 5.3 we would need the following result, which can be proved by analogous arguments.

\begin{lem}
Assume $\E \left[ | \log A |^3 \right] < \infty$, and set $\epsilon(m) = e^{(m-1)\mu(0)}$. Then there exist positive constants
$\widetilde{D}_1, \widetilde{c}, \widetilde{d}, {\mathfrak a}, {\mathfrak b}$
such that
\begin{equation}
\P \left\{ \epsilon(m) \le \Pi_{m-1} \le \widetilde{c} \, \epsilon(m), \quad \! Y_m \le \widetilde{d}, \quad \! {\mathfrak a} < B_1 < {\mathfrak b} \right\} \ge \frac{\widetilde{D}_1}{\sqrt{m}}
\end{equation}
for sufficiently large $m$.
\end{lem}

We now show by example the typical difficulty that one encounters when $B$ is allowed to be random.  In the proof of Lemma 5.3, we need to estimate
\[
\P \left\{ \gamma u \le \Pi_{k_u}, \quad \! \frac{\gamma b u}{i^2} e^l \le \Pi_{k_u - i} \max\left(1, B_{k(u)-i + 1}\right)
  \le \frac{\gamma b u}{i^2} \cdot e^{l+1} \right\}.
  \]
This estimate is obtained by considering
$\Pi_i^\prime := \Pi_{k_u}/\Pi_{k(u)-i}$,
and we need to have bounds on the two independent random variables, $\Pi_{k(u)-i}$ and $\Pi_i^\prime$.   For this purpose, we essentially need to eliminate the
$A_{k(u)-i+1}$ and $B_{k(u)-i+1}$ terms, estimating the above probability by
\begin{align*}
\sum_{j,r} \P \bigg\{ \gamma u & \le \Pi_{k_u-i} \Pi_{i-1}^\prime e^{j+1}, \quad \!  \frac{\gamma b u}{i^2} e^{-r-1} e^l \le \Pi_{k_u-i} \le  \frac{\gamma b u}{i^2} e^{-r} e^{l+1},\\[-.2cm]
   & \hspace*{1.5cm} e^j \le \max(1, A_{k_u - i + 1} ) \le e^{j+1}, \quad \! e^r \le \max(1, B_{k_u-i + 1}) \le e^{r+1} \bigg\},
\end{align*}
and then using that $\E \left[ A^\alpha + |B|^\alpha \right] < \infty$, to sum over all $j$ and $r$.  We omit the details, which are straightforward but technical.

%
%
%
%

\subsection{Proof of Theorem \ref{thml2}}
Let $b\geq 1$ and $B\in (0,b)$.
Then Theorem \ref{thml2} is a consequence of the following.

\begin{lem}
\label{lem:last}
Assume that the hypotheses of Theorem \ref{thml2} are satisfied.
Then there exists a constant $\theta \in (0,1)$ and ${\cal D}^\prime < \infty$ such that for every $\epsilon \in(0,1/2)$ and $u$ sufficiently large,
\begin{equation} \label{bbb1}
\P\left\{  Y_{n_u +k-1}\in\left( (1-\eps)u, \Big(1-\frac{\eps}{2} \Big)u \right),\quad \! \Pi_{n_u+k-1}> \frac{\eps}{2b} u
 \right\} \le  \eps^{1-\theta} \: \frac{{\cal D}^\prime \lambda^k (\alpha)}{\sqrt{\log u}} \: u^{-I(\tau)},
\end{equation}
where $n_u := \lfloor \tau \log u \rfloor$ and $k$ is any nonnegative integer, and the above result holds uniformly in $k$.
\end{lem}

Before presenting the proof of the lemma, we first show how it may be applied to establish the theorem.

\medskip

\noindent
{\bf Proof of Theorem \ref{thml2}.}
The proof is a simple consequence of the lemma.  Set $r_u = n_u +k-1$.  Since
\[
\big(0,u\big] = \bigcup_{i \ge 0}  \bigg( \bigg( 1 - \frac 1{2^i} \bigg) u, \bigg( 1 - \frac 1{2^{i+1}} \bigg) u  \bigg],
\]
and since $Y_{r_u+1} = Y_{r_u} + \Pi_{r_u} B_{r_u+1}$, it follows that
\begin{eqnarray*}
\P\left\{ Y_{r_u} \le u  \mbox{ and } Y_{r_u+1} > u \right\} &\le&
\sum_{i\ge 0}\P\left\{ Y_{r_u} \in \bigg( \bigg( 1 - \frac{1}{2^i} \bigg) u, \bigg( 1 - \frac{1}{2^{i+1}} \bigg) u  \bigg], \quad \! \Pi_{r_u}B_{r_u+1}>\frac u{2^{i+1}}  \right\} \\[.2cm]
&\le&
\sum_{i\ge 0}\P\left\{ Y_{r_u} \in \bigg( \bigg( 1 - \frac 1{2^i} \bigg) u, \bigg( 1 - \frac 1{2^{i+1}} \bigg) u  \bigg], \quad \! \Pi_{r_u}>\frac u{2^{i+1}b}  \right\} \\[.2cm]
&\le& \sum_{i\ge 0} \left( \frac{1}{2^{1-\theta}} \right)^i   \frac{{\cal D}^\prime \lambda^k (\alpha)}{\sqrt{\log u}} \: u^{-I(\tau)} \\[.2cm]
&\le& \frac{{\cal D} \lambda^k (\alpha)}{\sqrt{\log u}} \:u^{-I(\tau)},
\end{eqnarray*}
as required.
\halmos

\newcommand{\ej}{\frac{e^j}{10j^2}}

\medskip

\noindent
{\bf Proof of Lemma \ref{lem:last}.}
We begin by establishing the following result.

\begin{Assert}  For any $c>1$ and $\epsilon \in (0,1/2)$, there exist positive constants $\theta \in (0,1)$ and ${\cal D}^\prime < \infty$ such that,
for some finite constant $U$,
\begin{equation}
\label{eq: pt1}
\P\left\{ c^{-1}u \le Y_{n_{u+k-1}} \le cu,\quad \! \Pi_{n_u+k-1}> \eps u \right\}
\le \epsilon^{-\theta} \: \frac{{\cal D}^\prime \lambda^k (\alpha)}{\sqrt{\log u}} \: u^{-I(\tau)}, \quad \mbox{\rm all } k \ge0, \:\: u \ge U.
\end{equation}
\end{Assert}

\medskip

\noindent
{\bf Proof of the Assertion.}
Fix $k$, and set $r_u = n_u +k-1.$
Then define the set of indices
$$
{\cal W}_j^u = \Big\{ i:\; i< r_u  \quad \mbox{and} \quad (cu) e^{-j} \le \Pi_{i} B_{i+1} \le (cu) e^{-j+1}  \Big\}.
$$

Now suppose that $c^{-1}u \le Y_{r_u} \le cu$.  Then we claim that for some $j$, the number of elements in the set ${\cal W}_j^u$ must be greater
than $e^j/\left( 10 c^2 j^2 \right)$.  Indeed, if this were not the case, then (setting $\Pi_0 = 1$) we would have
\begin{eqnarray*}
Y_{r_u} &:=& \sum_{i=0}^{r_u-1}\Pi_i B_{i+1}\:\:\le \:\: \sum_{j=1}^\8\sum_{i\in {\cal W}_j^u}\Pi_i B_{i+1}\\
&& \hspace*{1cm} \mbox{(since }Y_{r_u} \le cu \Longrightarrow \Pi_i B_{i+1} \le cu \mbox{ for all }i \le r_u)\\[.2cm]
&\le& \sum_{j=1}^\8   \frac{e^j}{10c^2j^2} \cdot \frac{e c u}{e^j} \le \frac e{10 }\cdot \frac{\pi^2}{6}\cdot \frac{u}{c} < \frac{u}{c},
\end{eqnarray*}
a contradiction.

We now focus on the event $\{ c^{-1} u \le Y_{r_u} \le cu \}$, which appears on the left-hand side of \eqref{eq: pt1}.
Let ${\cal K}^u$ denote the following set of indices:
{ $$
{\cal K}^u = \bigg\{ (j,m_1,m_2):\ {j\ge 1},\ {1\le m_1 <r_u},\ {m_1+\frac{e^j}{10 c^2 j^2} < m_2 <r_u}.
\bigg\}
$$
Recall that for some $j$, ${\cal W}^u_j$ contains at least $e^j/\left( 10 c^2 j^2 \right)$ members.  This means that the first and last
occurrences  of the event described in ${\cal W}^u_j$ must be separated by a distance of at least $e^j/\left( 10 c^2 j^2 \right)$; that is, there must exist
values $m_1$ and $m_2$ such that
\[
(cu) e^{-j} \le \Pi_{m_i} B_{m_i+1} \le (cu) e^{-j+1}, \quad i=1,2, \qquad \mbox{and} \quad\!  m_2 - m_1 > \frac{e^j}{10 c^2 j^2}.
\]
Consequently,
\begin{align} \label{ee1}
\P\Big\{ & c^{-1}u  \le Y_{r_u} \le cu,\quad \! \Pi_{r_u}> \eps u \Big\} \nonumber\\[.1cm]
&\le \sum_{(j,m_1,m_2)\in {\cal K}^u}
\P\  \left\{(cu)e^{-j} \le  \Pi_{m_i} B_{m_i+1} \le (cu)e^{-j+1}, \:\: i=1,2; \:\: \Pi_{r_u}> \eps u  \right\} \nonumber\\[.2cm]
& \le \sum_{(j,m_1,m_2)\in {\cal K}^u}  \P\  \left\{ \Pi_{m_1}  \ge \frac{cu}{b} e^{-j} \right\}  \P \left\{ \Pi_{m_2 - m_1} B_1^{-1} \ge \frac{1}{b} e^{-1} \right\} \P \left\{ \Pi_{r_u - m_2} B_1^{-1} > \frac{\epsilon}{bc} e^{j-1} \right\}  \nonumber\\[.2cm]
&  := \sum_{(j,m_1,m_2)\in {\cal K}^u} P_1^u P_2^u P_3^u,
\end{align}
where $P_1^u, P_2^u, P_3^u$ denote, respectively, the three probabilities appearing in the previous expression on the right-hand side.

While we will ultimately need a sharper estimate, we first estimate these probabilities via Chebyshev's inequality by choosing parameters
$\beta_1 \in (0,\alpha)$ and $\beta_2 \in (0, \alpha \wedge 1)$ such that
\[
\rho_1 := \frac{\lambda(\beta_1)}{\lambda(\alpha)} < 1 \qquad \mbox{and} \qquad \rho_2 := \frac{\lambda(\beta_2)}{\lambda(\alpha)} < 1.
\]
Note that the parameter $\beta_2$ exists due to the assumption \eqref{eq: count2}.
Applying Chebyshev's inequality with the parameter $\alpha$ to the probability $P_1^u$, with parameter $\beta_1$ to  the probability $P_2^u$, and
with parameter $\beta_2$ to the probability $P_3^u$, we obtain by \eqref{ee1} that
\begin{align} \label{ee2}
\P\Big\{ c^{-1}u  &\le Y_{r_u} \le cu,\quad \! \Pi_{r_u}> \eps u \Big\} \nonumber\\[.1cm]
& \le C_1 \eps^{-\beta_2} \sum_{(j,m_1,m_2)\in {\cal K}^u} e^{j(\alpha - \beta_2)}  \Big( \lambda^{m_1}(\alpha) \lambda^{m_2-m_1}(\beta_1)
    \lambda^{r_u - m_2}(\beta_2) \Big) u^{-\alpha}  \nonumber\\[.2cm]
& = \Big( C_1 \eps^{-\beta_2} \Big) \lambda^k(\alpha) u^{-I(\tau)} \sum_{(j,m_1,m_2)\in {\cal K}^u}  e^{j(\alpha - \beta_2)} \rho_1^{m_2-m_1} \rho_2^{r_u-m_2}
\end{align}
for some constant $C_1 < \infty$, where we have used the assumption that $\lambda_B(-\alpha) < \infty$.

Next, fix $t> 0$ and divide the set ${\cal K}^u$ into four subsets, as follows:
\begin{eqnarray*}
{\cal K}^u_1 &=& \Big\{ (j,m_1,m_2)\in {\cal K}^u:\ e^j >  t r_u   \Big\};\\[.05cm]
{\cal K}^u_2 &=& \Big\{ (j,m_1,m_2)\in {\cal K}^u:\ e^j \le t r_u,\ m_2 <  r_u-r_u^{1/4}  \Big\};\\[.05cm]
{\cal K}^u_3 &=& \Big\{ (j,m_1,m_2)\in {\cal K}^u:\ e^j \le t r_u,\  m_1 <  r_u-2r_u^{1/4},\  m_2 \ge  r_u-r_u^{1/4}  \Big\};\\[.05cm]
{\cal K}^u_4 &=& \Big\{ (j,m_1,m_2)\in {\cal K}^u:\ e^j \le t r_u,\  m_1 \ge  r_u-2 r_u^{1/4}   \Big\}.
\end{eqnarray*}
We now study \eqref{ee2} by calculating the sum on the right-hand side separately over the respective sets ${\cal K}_i^u$, $i=1,\ldots,4$.

\medskip

\textbf{\textit{Case 1:}}  First, we estimate the sum over ${\cal K}_1^u$.
Since ${\cal K}^u_1 \subset {\cal K}^u$, we have $m_2-m_1 > e^j/\left(10c^2 j^2 \right) \ge L e^{j/2}$ for some constant $L>0$.
Thus for some
positive constant $L_1$,
\begin{align} \label{ee5}
\sum_{(j,m_1,m_2)\in {\cal K}^u_1} P_1^u P_2^u P_3^u
&\: \le \: \left( C_1 \eps^{-\beta_2} \right) \lambda^k(\alpha) u^{-I(\tau)}
   \sum_{(j,m_1,m_2)\in {\cal K}^u_1}  e^{j(\alpha - \beta_2)} e^{-L_1 e^{j}} \rho_2^{r_u-m_2} \nonumber\\[.2cm]
& \hspace*{2cm} =  o \left( \frac{\lambda^k(\alpha)}{\sqrt{\log u}} u^{-I(\tau)} \right)  \quad \mbox{as} \quad u \to \infty
\end{align}
when $t$ is sufficiently large.
The last step  follows since $\rho_2 < 1$ and, by definition, the set ${\cal K}^u$ contains at most $r_u := \lfloor \tau \log u \rfloor + k$ members,
while the subset ${\cal K}_u^1$ contains only those members where $e^j > t  r_u$ (so that in \eqref{ee1}, the sum over $j$ is finite and dominated by its initial term, that is,
$\sum_{{\cal K}_1^u} e^{j (\alpha - \beta_2)} e^{-L_1e^j} \le C_2 \big[ e^{j (\alpha - \beta_2)} e^{-L_1e^j} \big]_{\{ j=\log(tr_u)\}} \downarrow 0$ as $u \to \infty$).

\medskip

\textbf{\textit{Case 2:}}  Next consider the sum over ${\cal K}_2^u$.
In this case, $r_u - m_2 > r_u^{1/4}$ and so in \eqref{ee1},
\begin{align} \label{ee100}
\sum_{(j,m_1,m_2)\in {\cal K}^u_2} P_1^u P_2^u P_3^u
& \le \Big( C_1 \eps^{-\beta_2} \Big) \lambda^k(\alpha) u^{-I(\tau)} \sum_{(j,m_1,m_2)\in {\cal K}_2^u}  (t r_u)^{\alpha - \beta_2} \rho_1^{m_2-m_1} e^{-L_2 n^{1/4}}\nonumber\\[.2cm]
& \hspace*{2.5cm} =  o \left( \frac{\lambda^k(\alpha)}{\sqrt{\log u}} u^{-I(\tau)} \right)  \quad \mbox{as} \quad u \to \infty
\end{align}
for $L_2$ a positive constant and $\rho_1 < 1$, where we have again used that ${\cal K}^u$ contains  at most $r_u := \lfloor \tau \log u \rfloor + k$ members,
and on  the subset ${\cal K}^u_2$, we have $e^j \le t r_u$.

\medskip

\textbf{\textit{Case 3:}}  For the sum over ${\cal K}_3^u$, we can follow the same argument as in Case 2.   In \eqref{ee1}, we now utilize that $m_2 - m_1 > r_u^{1/4}$ and observe
that $\rho_2 < 1$ (rather than observing that $r_u - m_2 > r_u^{1/4}$ and $\rho_1 < 1$).  Hence, in either case, we have that
\[
\rho_1^{m_2 - m_1} \rho_2^{r_u - m_2} \le e^{-L_2 n^{1/4}},
\]
and \eqref{ee1} can be applied to deduce the same estimate as in \eqref{ee100}.

\medskip

\textbf{\textit{Case 4:}}  Finally, we estimate the sum over ${\cal K}_4^u$.  This estimate requires a more intricate calculation than \eqref{ee1}, relying now on Petrov's Theorem 3.1.

Since $m_1 \ge r_u - 2r_u^{1/4}$, we may apply Theorem  \ref{thm:3.1Petrov} to obtain that, uniformly in $m_1 \in [r_u - 2r_u^{1/4},r_u]$,
\begin{equation}
P_1^u :=  \P\  \left\{ \Pi_{m_1}  \ge \frac{cu}{b} e^{-j} \right\} \le \frac{C_3 e^{\alpha j}}{\sqrt{m_1}} \lambda^{m_1}(\alpha) u^{-\alpha}, \qquad u \ge U_0,
\end{equation}
independent of $k$,
where $C_3$ are $U_0$ are finite positive constants.
Thus, repeating the calculation in \eqref{ee2}, but using this estimate for $P_1^u$ in place of the previous estimate (which was based on Chebyshev's inequality),
we obtain that
\begin{eqnarray}
\sum_{(j,m_1,m_2)\in {\cal K}^u_4} P_1^u P_2^u P_3^u
\le \Big( C_4 \eps^{-\beta_2} \Big) \frac{\lambda^k(\alpha)}{\sqrt{\log u}} \: u^{-I(\tau)} \sum_{(j,m_1,m_2)\in {\cal K}^u_4}  e^{j(\alpha - \beta_2)} \rho_1^{m_2-m_1} \rho_2^{r_u-m_2}
\end{eqnarray}
for some finite constant $C_4$ and $u$ sufficiently large.  To complete the proof, it is sufficient to justify that the last sum is bounded.  For this purpose, first recall that
since ${\cal K}_4^u \subset {\cal K}^u$, then as argued in Case 1, we have that $m_2 - m_1 > L e^{j/2}$ for some $L>0$.  Hence for some positive constant $L_3$,
\begin{align}
\sum_{(j,m_1,m_2)\in {\cal K}^u_4} &  e^{j(\alpha - \beta_2)} \rho_1^{m_2-m_1} \rho_2^{r_u-m_2}
\le \sum_{(j,m_1,m_2)\in {\cal K}^u_4}  e^{j(\alpha - \beta_2)} e^{-L_1 e^{j/2}} \rho_1^{(m_2-m_1)/2} \rho_2^{r_u-m_2} \nonumber\\[.2cm]
& \le \left( \sum_j e^{j(\alpha - \beta_2)} e^{-L_1 e^{j/2}} \right) \left( \sum_{m_1 < m_2} \rho_1^{(m_2-m_1)/2} \right)
  \left(\sum_{m_2 < r_u}\rho_2^{r_u-m_2} \right) < \infty,
\end{align}
since $\rho_1 < 1$ and $\rho_2 < 1$.
Combining the estimates in Steps 1-4, we obtain \eqref{eq: pt1}, as required.
\halmos

\medskip

Returning now to the proof of the lemma,
set
\begin{eqnarray*}
J_{\eps} &=& \left( (1-\eps)u, \left(1-\frac{\eps}{2}\right)u \right), \quad \epsilon > 0;\\
Y_n^\prime &=& B_2 + \sum_{i=3}^{n+1} \left( A_2 \dots A_{i-1}\right) B_{i}, \quad n=1,2,\ldots;\\
\Pi_n^\prime &=& \prod_{i=2}^{n+1} A_i, \quad n=1,2,\ldots.
\end{eqnarray*}
Then for all $n$, $(Y_n, \Pi_n) \stackrel{\cal D}{=} (Y_n^\prime, \Pi_n^\prime)$ and
$Y_n = B_1 + A_1 Y_{n-1}^\prime$.

Suppose that the constant $a$ has been chosen such that w.p.1, the support of the law of $A_1$ is contained in the interval $[1/a,a]$.  Setting
$r_u = n_u + k -1$, we then obtain
\begin{align} \label{f0}
\P\left\{ Y_{r_u}\in J_{\eps}, \quad \! \Pi_{r_u}>\frac{\eps}{2b} u  \right\}
&\le \P\left\{  B_1+ A_1 Y^\prime_{r_u-1} \in J_{\eps}, \quad\! \Pi^\prime_{r_u-1} > \left(\frac{\eps}{2ab} \right) u  \right\} \nonumber\\[.1cm]
& \hspace*{-1cm} \leq \P\left\{ A_1\in \frac{1}{Y^\prime_{r_u-1}}\left( (1-\eps )u-b, \left(1-\frac{\eps}{2} \right)u\right),\    \Pi^\prime_{r_u-1}> \left(\frac{\eps}{2ab} \right) u \right\},
\end{align}
where $A_1$ is independent of $(Y^\prime_{r_u-1}, \Pi^\prime_{r_u-1})$. Moreover, since $a^{-1} \leq A_1\leq a$, we also have when $Y_{r_u}\in J_{\eps}$ that
\begin{equation} \label{f1}
Y^\prime_{r_u-1} \in \frac 1{A_1}\left( (1-\eps )u-b, \left(1-\frac{\eps}{2} \right)u\right) \subset \bigg( \frac{(1-\eps)u-b}{a},\: a\left(1-\frac{\eps}{2}\right) u  \bigg)\subset
\left(\frac{u}{2a}, au \right)
\end{equation}
for sufficiently large $u$, independent of $k$. Then
for fixed $Y^\prime_{r_u-1} \in \left(u/2a, au \right)$,
an easy calculation shows that
the length of the interval
\[
 \frac{1}{Y^\prime_{r_u-1}}\left( (1-\eps )u-b, \left(1-\frac{\eps}{2} \right)u\right)
\]
is bounded above by $d\epsilon$ for some positive constant $d$.  Hence, returning to \eqref{f0}, we obtain that
\begin{eqnarray} \label{f2}
\P\left\{ Y_{r_u}\in J_{\eps}, \quad \! \Pi_{r_u}>\frac{\eps}{2b} u  \right\}
 &\le& \int_{u/2a}^{au} \P\left\{ A_1\in \frac 1s \left( (1-\eps )u-b, \left(1-\frac{\eps}{2} \right)u\right) \right\} \nonumber\\[-.1cm]
 && \hspace*{3.5cm} \cdot \, \P\left\{ \Pi^\prime_{r_u-1} > \left( \frac{\eps}{2ab} \right) u,\quad \!Y^\prime_{r_u-1}\in ds\right\} \nonumber\\[.2cm]
&\le& d\eps \: \P\left\{ \frac {u}{2a} \le Y_{r_u-1}\le au,\quad \! \Pi_{r_u-1 }> \left( \frac{\eps}{2ab} \right) u \right\} \nonumber\\[.2cm]
&\le& d\eps \: \P \left\{ \frac {u}{a} \le Y_{r_u-1}\le au,\quad \! \Pi_{r_u-1 }>  \epsilon^\ast u \right\}
\end{eqnarray}
for certain positive constants $\epsilon^\ast$ and $u$.  Applying \eqref{eq: pt1} to the last quantity on the right-hand side yields \eqref{bbb1}, as required.
\halmos

\nocite{VP65,JCAV13,JCAV13b,SA00}

{\small
\bibliographystyle{chicago}
\bibliography{Bbl}

\begin{thebibliography}{}

\bibitem[\protect\citeauthoryear{Alsmeyer}{Alsmeyer}{2003}]{GA03}
Alsmeyer, G. (2003).
\newblock On the {H}arris recurrence of iterated random {L}ipschitz functions
  and related convergence rate results.
\newblock {\em J. Theoret. Probab.\/}~{\em {\bf 16}}, 217--247.

\bibitem[\protect\citeauthoryear{Alsmeyer and Iksanov}{Alsmeyer and
  Iksanov}{2009}]{AI2009}
Alsmeyer, G. and A.~Iksanov (2009).
\newblock A log-type moment result for perpetuities and its application to
  martingales in supercritical branching random walks.
\newblock {\em Electron. J. Probab.\/}~{\em {\bf 14}}, 289--312.

\bibitem[\protect\citeauthoryear{Alsmeyer and Mentemeier}{Alsmeyer and
  Mentemeier}{2012}]{AM2012}
Alsmeyer, G. and S.~Mentemeier (2012).
\newblock Tail behaviour of stationary solutions of random difference
  equations: the case of regular matrices.
\newblock {\em J. Difference Equ. Appl.\/}~{\em {\bf 18}}, 1305--1332.

\bibitem[\protect\citeauthoryear{Arfwedson}{Arfwedson}{1955}]{GA55}
Arfwedson, G. (1955).
\newblock Research in collective risk theory. {P}art {I}{I}.
\newblock {\em Skand. Aktuarietidskr.\/}, 53--100.

\bibitem[\protect\citeauthoryear{Asmussen}{Asmussen}{2000}]{SA00}
Asmussen, S. (2000).
\newblock {\em Ruin Probabilities}.
\newblock River Edge, NJ: World Scientific.

\bibitem[\protect\citeauthoryear{Asmussen and Sigman}{Asmussen and
  Sigman}{1996}]{SAKS96}
Asmussen, S. and K.~Sigman (1996).
\newblock Monotone stochastic recursions and their duals.
\newblock {\em Probab. Th. Eng. Inf. Sc.\/}~{\em {\bf 10}}, 1--20.

\bibitem[\protect\citeauthoryear{Billingsley}{Billingsley}{1986}]{PB86}
Billingsley, P. (1986).
\newblock {\em Probability and Measure\/} (2nd ed.).
\newblock New York: John Wiley and Sons.

\bibitem[\protect\citeauthoryear{Bollerslev}{Bollerslev}{1986}]{TB86}
Bollerslev, T. (1986).
\newblock Generalized autoregressive conditional heteroskedasticity.
\newblock {\em J.\ Econometrics\/}~{\em {\bf 31}}, 307--327.

\bibitem[\protect\citeauthoryear{Brofferio and Buraczewski}{Brofferio and
  Buraczewski}{2014}]{BB2014}
Brofferio, S. and D.~Buraczewski (2014).
\newblock On unbounded invariant measures of stochastic dynamical systems.
\newblock {\em {\rm To appear in} Ann. Probab\/}.

\bibitem[\protect\citeauthoryear{Buraczewski}{Buraczewski}{2009}]{bur2009}
Buraczewski, D. (2009).
\newblock On tails of fixed points of the smoothing transform in the boundary
  case.
\newblock {\em Stochastic Process. Appl.\/}~{\em {\bf 119}}, 3955--3961.

\bibitem[\protect\citeauthoryear{Buraczewski, Damek, Guivarc'h, Hulanicki, and
  Urban}{Buraczewski et~al.}{2009}]{DarEwa2009}
Buraczewski, D., E.~Damek, Y.~Guivarc'h, A.~Hulanicki, and R.~Urban (2009).
\newblock Tail-homogeneity of stationary measures for some multidimensional
  stochastic recursions.
\newblock {\em Probab. Theory Relat. Fields\/}~{\em {\bf 145}}, 385--420.

\bibitem[\protect\citeauthoryear{Buraczewski, Damek, Guivarc'h, and
  Mentemeier}{Buraczewski et~al.}{2014}]{BDGM2014}
Buraczewski, D., E.~Damek, Y.~Guivarc'h, and S.~Mentemeier (2014).
\newblock On multidimensional {M}andelbrot cascades.
\newblock {\em J. Difference Equ. Appl\/}~{\em 20(11)}, 1523--1567.

\bibitem[\protect\citeauthoryear{Buraczewski, Damek, Mikosch, and
  Zienkiewicz}{Buraczewski et~al.}{2013}]{BDMZ}
Buraczewski, D., E.~Damek, T.~Mikosch, and J.~Zienkiewicz (2013).
\newblock Large deviations for solutions to stochastic recurrence equations
  under {K}esten's condition.
\newblock {\em Ann. Probab.\/}~{\em {\bf 41}}, 2755--2790.

\bibitem[\protect\citeauthoryear{Buraczewski, Damek, and
  Zienkiewicz}{Buraczewski et~al.}{2014}]{BDZ}
Buraczewski, D., E.~Damek, and J.~Zienkiewicz (2014+).
\newblock Precise tail asymptotics of fixed points of the smoothing transform
  with general weights.
\newblock {\em {\rm To appear in} Bernoulli\/}.

\bibitem[\protect\citeauthoryear{Carmona, Petit, and Yor}{Carmona
  et~al.}{2001}]{PCFPMY01}
Carmona, P., F.~Petit, and M.~Yor (2001).
\newblock Exponential functionals of {L}\'{e}vy processes.
\newblock In O.~E. Barndorff-Nielsen, T.~Mikosch, and S.~I. Resnick (Eds.),
  {\em L\'{e}vy Processes: Theory and Applications}, pp.\  41--55. Boston:
  Birkh\"{a}user.

\bibitem[\protect\citeauthoryear{Collamore}{Collamore}{1998}]{JC98}
Collamore, J.~F. (1998).
\newblock First passage times of general sequences of random vectors: a large
  deviations approach.
\newblock {\em Stochastic Process. Appl.\/}~{\em {\bf 78}}, 97--130.

\bibitem[\protect\citeauthoryear{Collamore}{Collamore}{2009}]{JC09}
Collamore, J.~F. (2009).
\newblock Random recurrence equations and ruin in a {M}arkov-dependent
  stochastic economic environment.
\newblock {\em Ann. Appl. Probab.\/}~{\em {\bf 19}}, 1404--1458.

\bibitem[\protect\citeauthoryear{Collamore and Vidyashankar}{Collamore and
  Vidyashankar}{2013a}]{JCAV13b}
Collamore, J.~F. and A.~N. Vidyashankar (2013a).
\newblock Large deviation tail estimates and related limit laws for stochastic
  fixed point equations.
\newblock In G.~Alsmeyer and M.~L\"{o}we (Eds.), {\em Random Matrices and
  Iterated Random Functions}, pp.\  91--117. Berlin: Springer.

\bibitem[\protect\citeauthoryear{Collamore and Vidyashankar}{Collamore and
  Vidyashankar}{2013b}]{JCAV13}
Collamore, J.~F. and A.~N. Vidyashankar (2013b).
\newblock Tail estimates for stochastic fixed point equations via nonlinear
  renewal theory.
\newblock {\em Stochastic Process. Appl.\/}~{\em {\bf 123}}, 3378--3429.

\bibitem[\protect\citeauthoryear{Dembo and Zeitouni}{Dembo and
  Zeitouni}{1993}]{ADOZ93}
Dembo, A. and O.~Zeitouni (1993).
\newblock {\em Large Deviations Techniques and Applications}.
\newblock Boston: Jones and Bartlett.

\bibitem[\protect\citeauthoryear{Ellis}{Ellis}{1984}]{RE84}
Ellis, R. (1984).
\newblock Large deviations for a general class of random vectors.
\newblock {\em Ann. Probab.\/}~{\em \bf 12}, 1--12.

\bibitem[\protect\citeauthoryear{Engle}{Engle}{1982}]{RE82}
Engle, R.~F. (1982).
\newblock Autoregressive conditional heteroscedasticity with estimates of the
  variance of {U}nited {K}ingdom inflation.
\newblock {\em Econometrica\/}~{\em {\bf 50}}, 987--1007.

\bibitem[\protect\citeauthoryear{Enriquez, Sabot, and Zindy}{Enriquez
  et~al.}{2009}]{NECSOZ09}
Enriquez, N., C.~Sabot, and O.~Zindy (2009).
\newblock A probabilistic representation of constants in {K}esten's renewal
  theorem.
\newblock {\em Probab. Theory Relat. Fields\/}~{\em {\bf 144}}, 581--613.

\bibitem[\protect\citeauthoryear{Geman and Yor}{Geman and Yor}{1993}]{HGMY93}
Geman, H. and M.~Yor (1993).
\newblock Bessel processes, {A}sian options, and perpetuities.
\newblock {\em Math. Finance\/}~{\em \bf 3}, 349--375.

\bibitem[\protect\citeauthoryear{Goldie}{Goldie}{1991}]{CG91}
Goldie, C.~M. (1991).
\newblock Implicit renewal theory and tails of solutions of random equations.
\newblock {\em Ann. Appl. Probab.\/}~{\em {\bf 1}}, 126--166.

\bibitem[\protect\citeauthoryear{Guivarc'h}{Guivarc'h}{1990}]{gui1990}
Guivarc'h, Y. (1990).
\newblock Sur une extension de la notion de loi semi-stable.
\newblock {\em Ann. Inst. H. Poincar\'e Probab. Statist.\/}~{\em {\bf 26}},
  261--285.

\bibitem[\protect\citeauthoryear{Guivarc'h and Le~Page}{Guivarc'h and
  Le~Page}{2013a}]{GL2013}
Guivarc'h, Y. and {\'E}.~Le~Page (2013a).
\newblock Homogeneity at infinity of stationary solutions of multivariate
  affine stochastic recursions.
\newblock In G.~Alsmeyer and M.~L\"{o}we (Eds.), {\em Random Matrices and
  Iterated Random Functions}, pp.\  119--135. Berlin: Springer.

\bibitem[\protect\citeauthoryear{Guivarc'h and Le~Page}{Guivarc'h and
  Le~Page}{2013b}]{GL2014}
Guivarc'h, Y. and {\'E}.~Le~Page (2013b).
\newblock On the homogeneity at infinity of the stationary probability for an
  affine random walk.
\newblock {\em {\rm HAL archives-ouvertes ID: hal-00868944
  (https://hal.archives-ouvertes.fr/hal-00868944)}\/}.

\bibitem[\protect\citeauthoryear{Kesten}{Kesten}{1973}]{HK73}
Kesten, H. (1973).
\newblock Random difference equations and renewal theory for products of random
  matrices.
\newblock {\em Acta Math.\/}~{\em \bf 131}, 207--248.

\bibitem[\protect\citeauthoryear{Kl\"{u}ppelberg and
  Kostadinova}{Kl\"{u}ppelberg and Kostadinova}{2008}]{CKRK08}
Kl\"{u}ppelberg, C. and R.~Kostadinova (2008).
\newblock Integrated insurance risk models with exponential {L}\'{e}vy
  investment.
\newblock {\em Insurance Math. Econom.\/}~{\em {\bf 42}}, 560--577.

\bibitem[\protect\citeauthoryear{Kl\"{u}ppelberg and
  Pergamenchtchikov}{Kl\"{u}ppelberg and Pergamenchtchikov}{2004}]{CKSP04}
Kl\"{u}ppelberg, C. and S.~Pergamenchtchikov (2004).
\newblock The tail of the stationary distribution of a random coefficient
  {AR}(q) model.
\newblock {\em Ann. Appl. Probab.\/}~{\em {\bf 14}}, 971--1005.

\bibitem[\protect\citeauthoryear{Letac}{Letac}{1986}]{GL86}
Letac, G. (1986).
\newblock A contraction principle for certain {M}arkov chains and its
  applications. random matrices and their applications.
\newblock {\em Proceedings of AMS-IMS-SIAM Joint Summer Research Conference
  1984. Contemporary Mathematics\/}~{\em \bf 50}, 263--273.

\bibitem[\protect\citeauthoryear{Liu}{Liu}{2000}]{QL00}
Liu, Q. (2000).
\newblock On generalised multiplicative cascades.
\newblock {\em Stochastic Process. Appl.\/}~{\em {\bf 86}}, 263--286.

\bibitem[\protect\citeauthoryear{Mikosch}{Mikosch}{2003}]{TM03}
Mikosch, T. (2003).
\newblock Modeling dependence and tails of financial time series.
\newblock In B.~Finkenst\"{a}dt and H.~Rootz\'{e}n (Eds.), {\em Extreme Values
  in Finance, Telecommunications, and the Environment}, pp.\  185--286. Boca
  Raton: Chapman and Hall.

\bibitem[\protect\citeauthoryear{Mirek}{Mirek}{2011}]{MM11}
Mirek, M. (2011).
\newblock Heavy tail phenomenon and convergence to stable laws for iterated
  {L}ipschitz maps.
\newblock {\em Probab. Theory Relat. Fields\/}~{\em {\bf 151}}, 705--734.

\bibitem[\protect\citeauthoryear{Nummelin}{Nummelin}{1984}]{EN84}
Nummelin, E. (1984).
\newblock {\em General Irreducible Markov Chains and Non--negative Operators}.
\newblock Cambridge: Cambridge University Press.

\bibitem[\protect\citeauthoryear{Nummelin and Tuominen}{Nummelin and
  Tuominen}{1982}]{ENPT82}
Nummelin, E. and P.~Tuominen (1982).
\newblock Geometric of {H}arris recurrent {M}arkov chains with applications to
  renewal theory.
\newblock {\em Stoch. Process. Appl.\/}~{\em {\bf 12}}, 187--202.

\bibitem[\protect\citeauthoryear{Nyrhinen}{Nyrhinen}{2001}]{HN01}
Nyrhinen, H. (2001).
\newblock Finite and infinite time ruin probabilities in a stochastic economic
  environment.
\newblock {\em Stochastic Process. Appl.\/}~{\em \bf 92}, 265--285.

\bibitem[\protect\citeauthoryear{Paulsen}{Paulsen}{2002}]{JP02}
Paulsen, J. (2002).
\newblock On {C}ram\'{e}r-like asymptotics for risk processes with stochastic
  return on investments.
\newblock {\em Ann. Appl. Probab.\/}~{\em {\bf 12}}, 1247--1260.

\bibitem[\protect\citeauthoryear{Petrov}{Petrov}{1965}]{VP65}
Petrov, V.~V. (1965).
\newblock On the probabilities of large deviations for sums of independent
  random variables.
\newblock {\em Theory Probab. Appl.\/}~{\em \bf 10}, 287--298.

\bibitem[\protect\citeauthoryear{Petrov}{Petrov}{1995}]{VP95}
Petrov, V.~V. (1995).
\newblock {\em Limit Theorems of Probability Theory}, Volume~4 of {\em Oxford
  Studies in Probability}.
\newblock Oxford: Clarendon Press.

\bibitem[\protect\citeauthoryear{Rockafellar}{Rockafellar}{1970}]{RR70}
Rockafellar, R.~T. (1970).
\newblock {\em Convex Analysis}.
\newblock Princeton: Princeton Univ. Press.

\bibitem[\protect\citeauthoryear{Roitershtein}{Roitershtein}{2007}]{AR07}
Roitershtein, A. (2007).
\newblock One-dimensional linear recursions with {M}arkov-dependent
  coefficients.
\newblock {\em Ann. Appl. Probab.\/}~{\em {\bf 17}}, 572--608.

\bibitem[\protect\citeauthoryear{Siegmund}{Siegmund}{1976}]{DS76a}
Siegmund, D. (1976).
\newblock The equivalence of absorbing and reflecting barrier problems for
  stochastically monotone {M}arkov processes.
\newblock {\em Ann. Probab.\/}~{\em {\bf 4}}, 914--924.

\bibitem[\protect\citeauthoryear{Varadhan}{Varadhan}{1984}]{SV84}
Varadhan, S. R.~S. (1984).
\newblock {\em Large Deviations and Applications}.
\newblock Philadelphia: SIAM.

\bibitem[\protect\citeauthoryear{Vervaat}{Vervaat}{1979}]{WV79}
Vervaat, W. (1979).
\newblock On a stochastic difference equation and a representation of
  non-negative infinitely divisible random variables.
\newblock {\em Adv. Appl. Prob.\/}~{\em \bf 11}, 750--783.

\end{thebibliography}
}

\vspace{.3cm}
\noindent
{\small
{\sc J.\ F.\ Collamore},
Department of Mathematical Sciences,
University of Copenhagen,
Universitetsparken 5,\\
\hspace*{.4cm} DK-2100 Copenhagen \O,
Denmark.\\
\hspace*{.4cm} E-mail: collamore$@$math.ku.dk

\vspace{.25cm}
\noindent
{\small
{\sc D.\ Buraczewski,
E.\ Damek,
J.\ Zienkiewicz},
Instytut Matematyczny,
Uniwersytet Wroclawski,\\
\hspace*{.4cm} 50-384 Wroclaw,
pl.\ Grunwaldzki 2/4,
Poland.\\
\hspace*{.4cm} E-mail:  dbura$@$math.uni.wroc.pl,
edamek$@$math.uni.wroc.pl,
zenek$@$math.uni.wroc.pl

\end{document}